\documentclass[a4paper,fleqn]{cas-dc}

\usepackage{amssymb}
\usepackage[numbers]{natbib}
\usepackage{hhline}
\usepackage{threeparttable}
\usepackage{tabularx}
\usepackage{multirow}
\usepackage{xcolor}
\usepackage{booktabs}
\usepackage{multicol}
\usepackage[ruled,vlined]{algorithm2e}
\usepackage[figuresright]{rotating} 
\usepackage{subfigure}
\usepackage{graphicx}
\usepackage{colortbl}
\usepackage{color}
\usepackage[numbered]{bookmark}
\usepackage{hyperref}
\hypersetup{
        bookmarksnumbered=true,
        bookmarksopen=true,
        bookmarksopenlevel=5,
        pdfstartview=Fit,
        pdfpagemode=UseOutlines
    }
\usepackage{array}
\usepackage[figuresright]{rotating}
\usepackage{threeparttable}

\usepackage[switch]{lineno}

\def\tsc#1{\csdef{#1}{\textsc{\lowercase{#1}}\xspace}}
\tsc{WGM}
\tsc{QE}
\tsc{EP}
\tsc{PMS}
\tsc{BEC}
\tsc{DE}

\begin{document}
\let\WriteBookmarks\relax
\def\floatpagepagefraction{1}
\def\textpagefraction{.001}
\shorttitle{FFEEA for Communication Satellite Network}
\shortauthors{Xuemei Jiang et~al.}

\title [mode = title]{An Evolutionary Task Scheduling Algorithm Using Fuzzy Fitness Evaluation Method for Communication Satellite Network}

\author[1]{Xuemei Jiang}
\ead{xuemeijane@126.com}
\address[1]{School of Electronic Engineering, Xidian University, Xi’an, China}
\author[2]{Yangyang Guo}
\ead{g2002dmu@163.com}
\address[2]{School of Systems Science, Beijing Jiaotong University, Beijing, China}
\author[3]{Yue Zhang}
\ead{zhangyue1127@buaa.edu.cn}
\address[3]{School of Reliability and Systems Engineering, Beihang University,Beijing, China}
%
\author[4]{Yanjie Song}[orcid=0000-0002-4313-8312]
\ead{songyj_2017@163.com}
\address[4]{National Engineering Research Center of Maritime Navigation System, Dalian Maritime University, Dalian, China}
%
\author[5,6,7]{Witold Pedrycz}
\ead{wpedrycz@ualberta.ca}
\address[5]{Department of Electrical and Computer Engineering, University of Alberta, Edmonton, Canada}
\address[6]{Systems Research Institute, Polish Academy of Sciences, Poland}
\address[7]{Faculty of Engineering and Natural Sciences, Department of Computer Engineering, Sariyer/Istanbul, Turkiye}
\author[1]{Lining Xing}
\ead{lnxing@xidian.edu.cn}


\begin{abstract}
Communications satellite network (CSN), as an integral component of the next generation of communication systems, has the capability to offer services globally. Data transmission in this network primarily relies on two modes: inter-satellite communication and satellite-to-ground station communication. The latter directly impacts the successful reception of data by users. However, due to resource and task limitations, finding a satisfactory solution poses a significant challenge. The communication satellite-ground station network scheduling problem (CS-GSNSP) aims to optimize CSN effectiveness by devising a plan that maximizes link construction time while considering constraints associated with satellite operation modes. The large number of tasks and numerous constraints in the problem result in a time-consuming evaluation of fitness function values. To address this issue, we propose a fuzzy fitness evaluation method (FFEM) that employs fuzzy or real evaluation methods based on individual similarity degrees. \textcolor[rgb]{0,0,0}{Additionally, we introduce an evolutionary algorithm based on FFEM, called evolutionary algorithm based on FFEM (FFEEA), for iteratively searching high-quality network construction schemes.} In FFEEA, an adaptive crossover approach is used for efficient population search. Finally, extensive experiments are conducted to demonstrate that our proposed fuzzy fitness evaluation method and other improvement strategies significantly enhance satellite network service time. \textcolor[rgb]{0,0,0}{The study introduces a novel approach to enhance the efficiency of solving combinatorial optimization problems, such as CS-GSNSP, by mitigating the complexity associated with fitness evaluation.}

\end{abstract}




\maketitle

\section{Introduction}

In recent years, communication satellites have exhibited a more rapid development pace than anticipated by the general public and have emerged as a novel paradigm of network communication. The communication satellite network (CSN) established by SpaceX enables global connectivity and facilitates seamless integration into our daily lives \cite{spangelo2015optimization}. Consequently, numerous new industries reliant on CSN for product or service provision have also surfaced. The increasing demand for diverse communication applications, encompassing CSN and its associated sectors, necessitates stringent requirements for the construction of communication satellite systems.

The communication satellite system is responsible for managing the normal on-orbit operation of communication satellites and assigning them to complete various communication tasks according to mission requirements \cite{lee2008task}. Within this system, the execution of communication tasks relies on establishing effective links between satellites and ground stations. The success of data upload (or transmission) depends directly on the quality of these links. However, due to fixed orbits, limited ground station locations, equipment availability, and inherent capabilities, there are constraints on the total duration available each day for constructing communication links \cite{khojah2022multi}. Communication satellites have two operational modes to alleviate resource limitations. In addition to the regular mode, a new feed switching (FS) operation mode allows immediate task execution after completing a task. Nevertheless, certain tasks may still remain unfulfilled due to various reasons.\textcolor[rgb]{0,0,0}{ Therefore, in order to optimize communication satellite-ground station networking scheduling problem (CS-GSNSP), it is crucial to maximize the overall duration of constructed communication links by equationting an efficient scheme based on satellite and ground station resources.}

\textcolor[rgb]{0,0,0}{Currently, there is limited research on the CS-GSNSP, whereas the satellite range scheduling problem (SRSP), a well-studied classic in satellite scheduling for the past two decades, exhibits numerous similarities. Both problems involve establishing antenna capture and link construction between satellites and ground stations. However, the SRSP problem encompasses a broader scope as it includes communication satellites, Earth observation satellites (EOSs), and navigation satellites that require telemetry, tracking, and commanding (TT\&C).} Chen et al. \cite{chen2021population} and Xiang et al. \cite{xiang2024hierarchical} have developed mixed integer programming models to maximize task profits. Song et al. \cite{song2019learning} also proposed a dual-objective optimization model considering additional objectives of minimizing task failure rate. Furthermore, Du et al. \cite{du2019moea} have defined a specific objective and proposed a three-objective optimization model. The CS-GSNSP mathematical planning model will be further designed by combining the requirements of constructing links with the characteristics of communication satellites, based on the existing SRSP problem-related models. In this particular scenario, where communication satellite antennas possess stronger capabilities allowing for quick establishment of communication links with new ground station antennas using FS operating mode, this additional work mode increases the total time available for building communication links but also introduces complexity in resource scheduling. Consequently, algorithm design becomes highly demanding necessitating the selection of solutions capable of long-term link establishment along with task execution plans from numerous task combination relations.

The complexity of the satellite scheduling problem is extremely high. Both the Earth observation satellite scheduling problem (EOSSP) and SRSP are classified as NP-Hard. When using a precise solution algorithm, it inevitably faces the challenge of exponential growth in computing time. For large-scale satellite resource scheduling problems, employing a precise solution algorithm is not suitable. To successfully address various satellite scheduling problems such as EOSSP and SRSP, researchers have proposed numerous heuristic algorithms and meta-heuristic algorithms. Among these algorithms, genetic algorithm, ant colony algorithm, particle swarm algorithm, and other meta-heuristic algorithms have demonstrated significant advantages through population search techniques. Although meta-heuristic algorithms can handle large-scale optimization problems, they still face poor search effectiveness even after extensive search iterations. Several improvement strategies have been employed to enhance the search performance of these algorithms. However, their search mechanisms cannot guarantee optimal solutions. When applying a meta-heuristic algorithm to solve complex engineering application problems like CS-GSNSP, finding an optimal solution becomes nearly impossible due to multiple constraints that result in high computational costs. Therefore, the primary focus lies in quickly obtaining high-quality solutions with minimal loss during practical applications. Throughout the entire search process of a meta-heuristic algorithm, fitness evaluation (FE) constitutes the most computationally expensive component. In this study, we propose a fuzzy fitness evaluation method (FFEM) by incorporating fuzzy concepts into FE methodology for solving CS-GSNSP using an evolutionary algorithm (EA), named evolutionary algorithm based on FFEM (FFEEA). FFEEA adaptively selects fitness evaluation and population evolution methods based on specific problem characteristics. \textcolor[rgb]{0,0,0}{The FFEM used in FFEEA can effectively reduce the cost of population fitness evaluation and improve the search efficiency.} The primary original contributions of this paper are as follows.

1. A mixed integer programming model is constructed to address the CS-GSNSP, taking into consideration the intricate operational scenarios of communication satellites in both regular and feed switching modes. The objective function is established to maximize the profit of link construction, aligning with the construction objectives of a communication satellite system. Furthermore, constraints are incorporated in the model to account for factors such as link construction task execution requirements and working capabilities of communication satellite and ground station equipment.

2. We propose a genetic algorithm that uses a fuzzy method to evaluate fitness function values. The fuzzy method enables quick fitness evaluation based on the membership relationship between individuals and optimal individuals. Our algorithm incorporates an $\varepsilon$-evaluation strategy selection mechanism, allowing individuals to choose between actual and fuzzy fitness evaluations. To dynamically determine the search mode used in the next-generation population search, we employ an adaptive crossover method that integrates multiple crossover operators based on search performance. Additionally, our algorithm adopts an elite retention strategy to accelerate convergence.

3. The experiment validates the efficacy of the proposed FFEEA in addressing the CS-GSNSP from multiple perspectives. \textcolor[rgb]{0,0,0}{In comparison to state-of-the-art algorithms, FFEEA demonstrates accelerated population iterative search through its novel fitness evaluation method.} This fitness evaluation approach not only proves suitable for the problem under investigation in this study but also exhibits potential applicability to other combinatorial optimization problems.

\textcolor[rgb]{0,0,0}{The rest of this article is organized as follows. The following section describes related research. Section \ref{Model} introduces the problem statement of CS-GSNSP and provides the corresponding mathematical programming model. Section \ref{Proposed Method} introduces the overall framework of the genetic algorithm based on the fuzzy fitness evaluation method, along with its main evolutionary strategies. In the section \ref{Experimental Studies}, an extensive set of experiments is designed to validate the effectiveness of the proposed algorithm. Section \ref{Conclusion}  analyze findings and discuss future research directions.}

\section{Related Work}
\label{Related Work}
This section presents the current research status of the communication satellite scheduling problem and explores the application of evolutionary algorithms for solving such problems.

\subsection{Communication Satellite Scheduling Problem}

The communication satellite scheduling problem (CSSP) is a complex combinatorial optimization problem, similar to the classical problems of vehicle routing problem (VRP) and traveling salesman problem (TSP) \cite{barbulescu2004scheduling}. Research on these types of problems plays a crucial role in improving the operational efficiency of communication satellites. Al-hraishawi et al. \cite{al2021scheduling} focused on generating high-quality communication satellite scheduling schemes through carrier aggregation and proposed a novel load-balancing scheduling algorithm. Peng et al. \cite{peng2021hybrid} addressed CSSP from the perspective of satisfying user needs and optimizing resources by proposing a feasible optimization method based on the minimum mean square error (MMSE) criterion and logarithmic linearization to search for reasonable power allocation for user terminals. Wei et al. \cite{wei2008genetic} constructed models and designed algorithms for communication satellite data relay transmission scheduling problems while fully considering resource capacity constraints. \textcolor[rgb]{0,0,0}{Wang et al. \cite{wang2023satellite} proposed a dynamic spatio-temporal approximation (DSTA) model and beam collaboration algorithm for CSSP scheduling, which was comprehensively tested through experiments to verify its correctness and effectiveness.} Dai et al. \cite{dai2018satellite} aimed at optimizing terrestrial cellular network coverage using communication satellites by reducing the number of required satellites while effectively covering regions, whereas Liu et al. \cite{liu2022data} proposed an improved strategy combining rules with operators that adopted parallel adaptive large neighborhood search algorithms for communication satellite constellation routing.

\textcolor[rgb]{0,0,0}{Currently, the majority of research in communication satellite scheduling primarily focuses on user scheduling and hardware-level design of satellite communication systems (e.g., beamforming and bandwidth allocation).} Moreover, these studies emphasize the provision of reliable services to specific regions using communication satellites, while paying less attention to CSNs that resemble the Internet. However, as CSNs are expected to become the future mainstream global network, there exist numerous challenges impacting service quality that necessitate resolution through model building and algorithm design for addressing the CS-GSNSP.

\subsection{\textcolor[rgb]{0,0,0}{Evolutionary Algorithm to Solve Satellite Scheduling Problem}}

The evolutionary algorithm (EA) exhibits excellent performance in addressing various satellite scheduling problems. This section primarily examines the relevant literature on using EA to tackle communication satellite scheduling problems, with a particular focus on SRSP, which is closely aligned with this research.

Zhang and Xing \cite{zhang2022improved} proposed an improved genetic algorithm (IGA) to solve the satellite data down transmission task scheduling problem. IGA employed a novel approach to encode and decode both satellite and terrestrial resources. To enhance the algorithm's search efficiency, a set of heuristic operators was incorporated into the algorithm framework. Song et al. \cite{song2020knowledge} introduced a knowledge-based genetic algorithm for addressing relay satellite scheduling problems. The proposed algorithm adopted a knowledge-driven population initialization and change strategy. Experimental validation was conducted to assess the overall performance of the algorithm along with its improvement strategies. Jing et al. \cite{jing2023energy} put forward an energy-efficient genetic algorithm for laser communication routing in communication satellites, which significantly improved the service life of satellites without compromising network performance. In order to enhance the utilization efficiency of communication satellites, Liu et al. \cite{liu2022communication} presented an improved whale optimization algorithm (IWOA). A range of search and detection operators were employed within the IWOA framework, demonstrating that it can effectively improve user services in communication satellite systems.

\textcolor[rgb]{0,0,0}{The study of CS-GSNSP benefits from the reference significance of SRSP, which has witnessed successful applications of EA. Du et al. \cite{du2019moea} and Song et al. \cite{song2019learning} respectively proposed MOEA approaches for solving multi-objective (MO) SRSP problems.} However, these MOEAs incur higher computational costs compared to single objective (SO) EAs and pose challenges in their application to real-world task scheduling systems. Ref. \cite{song2024energy} considered the energy consumption factor in SRSP. Song et al. \cite{song2019learning} integrated a reinforcement learning-based neighborhood search strategy into the memetic algorithm (MA), incorporating various random, heuristic, and energy-efficient optimization techniques for problem design purposes. Simulation experiments demonstrated that the proposed algorithm effectively reduced power consumption during satellite telemetry, track, and command (TT\&C) processes. Chen et al. \cite{chen2021population} introduced a perturbation and elimination strategy based on genetic algorithms to address multi-star-multi-station SRSP problems while Ref.\cite{song2023cluster} combined clustering methods with genetic algorithm to tackle SRSP issues by designing clustering-based crossover and mutation operators that dynamically selected strategies according to the search process.

The aforementioned studies have demonstrated significant advancements in algorithm search processes and strategies across various aspects. These improvements primarily focus on detailed problem-specific designs. To the best of our knowledge, no prior researchers have employed the fuzzy fitness evaluation method within the domain of satellite scheduling. This evaluation approach not only reduces computational costs in CS-GSNSP but also offers potential solutions for other complex combinatorial optimization problems.

\section{Model}
\label{Model}

This section provides a comprehensive description of the CS-GSNSP, encompassing the symbols and variables employed in the model, assumptions, objective functions, and constraints.

\subsection{Problem Description}

\textcolor[rgb]{0,0,0}{In the CS-GSNSP, a set of communication satellite and ground station resources are used to construct the communication satellite network (CSN) \cite{kilic2019modeling}.} Each communication satellite can establish a communication link with a ground station antenna when it passes over it, using its onboard antenna. This period during which a communication link can be established is referred to as the visible time window (VTW). The satellite can only transmit data when there is an available VTW. Due to the fixed orbit of the satellite, the number of orbits around Earth within the planned time frame is limited, resulting in restricted availability of VTWs. The key to solving this problem is how to choose VTW for establishing communication links. The communication satellite can also perform tasks in either the regular working mode or the feed switching working mode based on the overlapping relationship between VTWs.

\begin{figure}[htp]
\centering
\subfigure[Schematic Diagram of Non-fulfillment of Time Requirements]{
\includegraphics[width=0.5\textwidth]{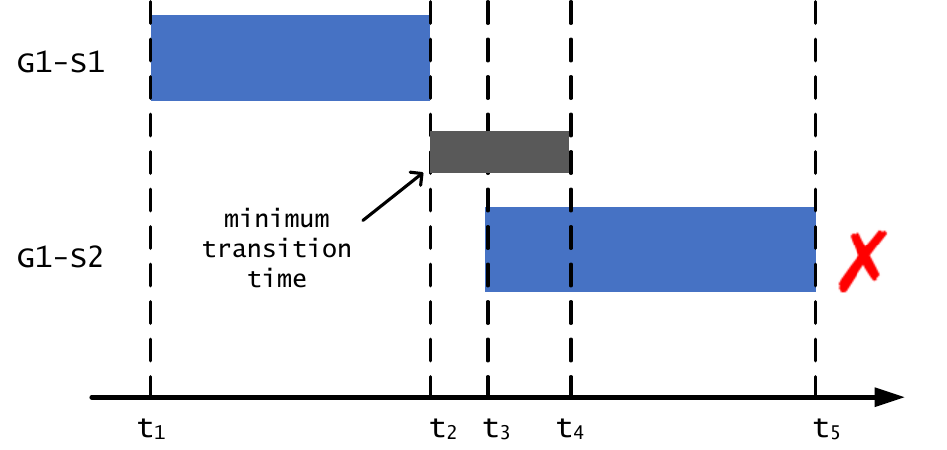}}
\subfigure[Schematic for Meeting Time Requirements]{
\includegraphics[width=0.5\textwidth]{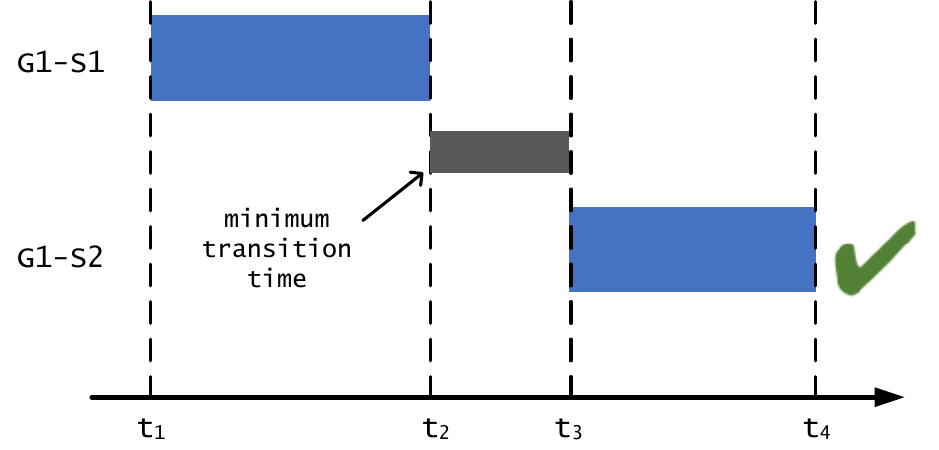}}
\caption{\textcolor[rgb]{0,0,0}{Schematic Diagram of Feeder Transition Time Requirements}}
\label{Schematic Diagram of Feeder Transition Time Requirements}
\end{figure}

Once a time window is selected, the specific execution time of the task needs to be determined based on the equipment capabilities of both the satellite and ground station. To ensure link stability, a certain setup time is required prior to establishing communication, as well as a certain processing time upon termination. A ground station antenna can only establish communication with a corresponding satellite antenna, and after completing a task, it must undergo equipment parameter switching for a designated period before commencing the next task. Figure \ref{Schematic Diagram of Feeder Transition Time Requirements}-(a) and Figure \ref{Schematic Diagram of Feeder Transition Time Requirements}-(b) depict schematic diagrams illustrating task switching at the ground station. In Figure \ref{Schematic Diagram of Feeder Transition Time Requirements}-(a), two tasks fail to meet the required transition time and thus cannot be executed. However, in Figure \ref{Schematic Diagram of Feeder Transition Time Requirements}-(b), both tasks successfully meet the ground station's transition time requirements.

\begin{figure}[htp]
\centering
\includegraphics[width=0.4\textwidth]{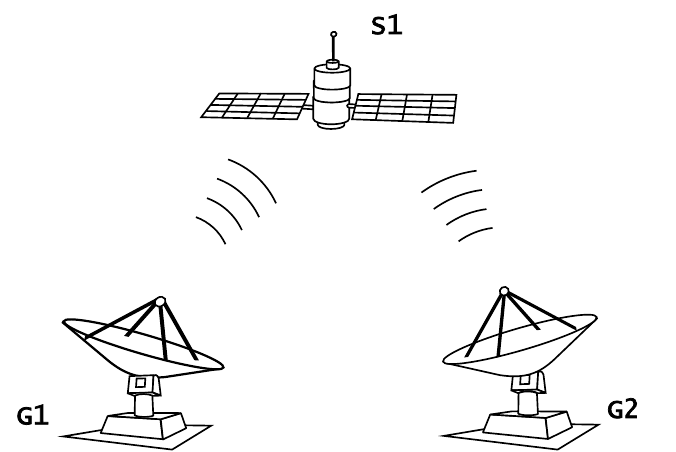}
\caption{Switching Diagram of the Ground Station Performing Two Tasks}
\label{Switching Diagram of the Ground Station Performing Two Tasks}
\end{figure}

\begin{figure}[htp]
\centering
\includegraphics[width=0.4\textwidth]{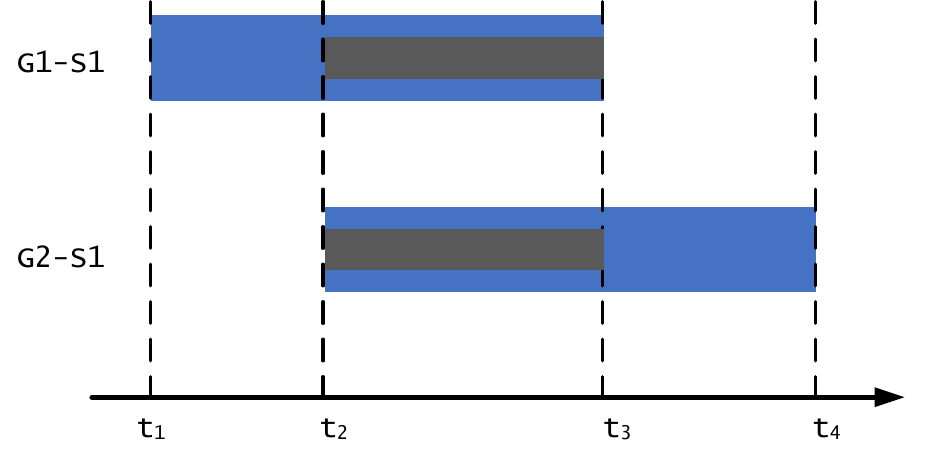}
\caption{Satellite Feed Switching Diagram}
\label{Satellite Feed Switching Diagram}
\end{figure}

A communication satellite is equipped with two antennas for establishing a communication link, but only one antenna is required to establish a connection with a ground station (refer to Figure \ref{Switching Diagram of the Ground Station Performing Two Tasks}). In the regular working mode, after completing a task, the satellite needs to configure several parameters before proceeding to the next task. However, by employing the feed switching operation mode, it becomes possible to directly switch from one ground station's link-building task to another without any waiting time. To ensure continuous connectivity, both tasks using this mode should have sufficient overlap time (typically less than the minimum time required for ground station task switching). Figure \ref{Satellite Feed Switching Diagram} illustrates a schematic diagram of satellite feed switching. This working mode significantly reduces wasted time during task transitions. Overall, addressing the CS-GSNSP involves developing an optimized work plan that adheres to constraints imposed by communication satellites and ground stations in order to achieve maximum CSN efficiency.

\subsection{Symbols and Variables}

In this section, we collect all symbols used in this study.

$S$: the set of communication satellites, the number of satellites is $|S|$, for satellite $sat_{i}$ contains the following attributes.

$A_i^s$: antenna set of satellite $sat_{i}$. The number of antennas is $\left|A_i\right|$. For antenna $a_{i}$, the following attributes are included.

-$\tau_{im}$: attitude adjustment time of an antenna.

-$\Delta_{imm^{\prime}}$: indicates the minimum overlap time for performing feed switching.

\textcolor[rgb]{0,0,0}{$G$: the set of ground stations, the number of ground stations is $|G|$, for the ground station $g_j$ contains the following attributes.}

-$A_{j}^{g}$: indicates the antenna set belonging to the ground station $g_{j}$. The number of antennas is $\left|A_j^g\right|$.

$\varphi_{ii^\prime}^j$: the minimum interval between tasks of ground station antenna switching (different star conditions at the same station).

$TW_{mn}^{ij}$: \textcolor[rgb]{0,0,0}{a set of time windows composed of satellite antenna $sat_{i}$ and ground station antenna $g_{j}$.} For the time window $tw_{k}$ in $TW_{mn}^{ij}$ contains the following attributes.

-$e\nu t_{k}$: the earliest visible time of the time window.

-$l\nu t_{k}$: the latest visible time of the time window.

$T_{i}$: a set of tasks with no time window overlap between satellite $sat_{i}$ and the ground station.

-$t_{o}$: link construction task between antenna $a_{i}$ and antenna $g_{j}$ in time window $tw_{k}$.

-$p_{o}$: unit time profit of link construction task $task_k^{o}$.

$\overline{T}_i$: the set of tasks with time window overlap between satellite $sat_{i}$ and the ground station.

-$\overline{t}_o$: construction task of feed switching link between antennas $a_{i}$ and $a_i^{\prime}$ and ground station antennas $g_{j}$ and $g_j^{\prime}$ in time window $tw_{k}$ and $tw_{k^{\prime}}$.

-$\overline{p}_o$: the unit time profit of the switchover link construction task $o$.

$M$: \textcolor[rgb]{0,0,0}{a large integer.}

\textbf{Decision variables:}

$x_{ijk}^{mno}$: whether satellite $sat_{i}$ and ground station $g_{j}$ are performing the regular link building task $t_{o}$ in the time window $tw_{k}$. If executed, $x_{ijk}^{mno}$=1. Otherwise, $x_{ijk}^{mno}$=0.

$st_{o}$: the start time of task $o$.

$et_{o}$: the completion time of task $o$.

$\overline{x_{ijj^{\prime}k}^{mno}}$: whether satellite $sat_{i}$ and ground stations $g_{j}$ and $g_j^{\prime}$ perform the feed switching link construction task $t_{o}$. If executed, $\overline{x_{ijj^{\prime}k}^{mno}}$=1. Otherwise, $\overline{x_{ijj^{\prime}k}^{mno}}$=0.

$\overline{st_o}$: the start time of task $o$.

$\overline{et_o}$: the completion time of task $o$.

\subsection{Mathematical Model}

\label{Mathematical Model}

\textcolor[rgb]{0,0,0}{Several specific assumptions for CS-GSNSP are provided as follows. The mathematical programming model will be constructed based on the following assumptions.}

\textcolor[rgb]{0,0,0}{\textbf{Assumptions:}}

1. All link construction tasks have been predetermined prior to task planning, ensuring that no individual tasks will be compelled to cancel due to external factors during the planning process.

2. Communication satellites and ground stations are capable of maintaining uninterrupted operational conditions at all times.

3. Each communication chain-building task is limited to a single execution, with no provision for multiple or periodic executions.

4. \textcolor[rgb]{0,0,0}{The satellites possess sufficient power reserves to successfully accomplish the link construction tasks.}

5. Once each link construction task has been effectively planned, it must be executed without fail.

In CS-GSNSP, our objective is to identify a task execution scheme that maximizes the profits derived from link construction activities, encompassing both regular and feed switching modes. The pursuit of higher profits implies a greater capacity to fulfill users' needs.

\textcolor[rgb]{0,0,0}{\textbf{Objective function:}}
\begin{equation}
\label{actual profit}
\max f_{reg}+f_{fee}
\end{equation}
\begin{equation}
 f_{reg}=\sum_{i\in S}\sum_{j\in G}\sum_{k\in TW_{mn}^{ij}}\sum_{m\in A_i^s}\sum_{n\in A_j^g}\sum_{o\in T_i}\left(et_o-st_o\right)\cdotp p_o\cdot x_{ijk}^{mno}
\end{equation}
\begin{equation}
 f_{fee}=\sum_{i\in S}\sum_{j\in G}\sum_{j'\in G}\sum_{k\in TW_{mn}^{ij}}\sum_{m\in A_i^s}\sum_{n\in A_j^s}\sum_{o\in T_i}\left(\overline{et_o}-\overline{st_o}\right)\cdotp\overline{p_o}\cdot\overline{x_{ijj^{\prime}k}^{mno}}
\end{equation}
where $f_{reg}$ represents the profits achievable in the regular operational mode, $f_{fee}$ denotes the profits attainable through switching to a feed-based working mode, $x_{ijk}^{mno}$ indicates whether the link construction task is executed in the regular operational mode, $\overline{x_{ijj^{\prime}k}^{mno}}$ determines whether to execute the link building task using feed switching mode.

\textcolor[rgb]{0,0,0}{\textbf{Constraints:}}

The communication satellite can adopt both the regular mode and the feed switching mode, which leads to variations in the restrictive conditions required for each mode to perform tasks. Therefore, this section will provide a detailed introduction to the general constraints and specific constraints associated with each working mode.

\textcolor[rgb]{0,0,0}{\textbf{-General constraints:}}

$\bullet$ The start time is earlier than the finish time.

\begin{equation}
    \begin{array}{l}
 st_o\cdot x_{ijk}^{mmo}\leq et_o,\forall i\in S,j\in G,\\
 k\in TW_{mm}^{ij},m\in A_i^s,n\in A_j^g,o\in T_i
    \end{array}
\end{equation}

\begin{equation}
\begin{array}{l}
 \overline{st_o}\cdot\overline{x_{ijj'k}^{mno}}\leq\overline{et}_o,\forall i\in S,j,j'\in G,\\
 k\in TW_{mn}^{ij},m\in A_i^s,n\in A_j^\mathrm{g},o\in\overline{T_i}
 \end{array}
\end{equation}

$\bullet$ The commencement of every satellite task is required to fall within the visible time window.
\begin{equation}
\begin{array}{l}
 e\nu t_k\cdot x_{ijk}^{mno}\leq st_o,\forall i\in S,j\in G,\\TT
 k\in TW_{mn}^{ij},m\in A_i^s,n\in A_j^g,o\in T_i
 \end{array}
\end{equation}
\begin{equation}
\begin{array}{l}
 evt_k\cdot\overline{x_{ijj'k}^{nmo}}\leq\overline{st_o},\forall i\in S,j,j'\in G,\\
 k\in TW_{mn}^{ij},m\in A_i^s,n\in A_j^\mathrm{g},o\in\overline{T_i}
 \end{array}
\end{equation}

$\bullet$ The completion of each satellite task is required to fall within the visible time window.
\begin{equation}
\begin{array}{l}
 et_o\cdot x_{ijk}^{mno}\leq l\nu t_k,\forall i\in S,j\in G,\\
 k\in TW_{mn}^{ij},m\in A_i^s,n\in A_j^g,o\in T_i
\end{array}
\end{equation}
\begin{equation}
\begin{array}{l}
 \overline{et}_o\cdot\overline{x_{ijjk}^{mno}}\leq l\nu t_k,\forall i\in S,j,j'\in G,\\
 k\in TW_{mn}^{ij},m\in A_i^s,n\in A_j^\mathrm{g},o\in\overline{T}_i
 \end{array}
\end{equation}

$\bullet$ The establishment of a communication link for a satellite necessitates the use of a single antenna at any given time.
\begin{equation}
\begin{array}{l}
 \sum_{m\in A_i^s}x_{ijk}^{mno}\leq1,\forall i\in S,j\in G,\\
 k\in TW_{mn}^{ij},n\in A_j^g,o\in T_i
 \end{array}
\end{equation}
\begin{equation}
\begin{array}{l}
 \sum_{m\in A_i^s}\overline{x_{ijj^{\prime}k}^{mno}}\leq1,\forall i\in S,j,j^{\prime}\in G,\\
 k\in TW_{mn}^{ij},n\in A_j^g,o\in\overline{T_i}
 \end{array}
\end{equation}

$\bullet$ The ground station antenna is capable of establishing a communication link with only one satellite antenna at any given time.
\begin{equation}
\begin{array}{l}
 \sum_{n\in A_j^g}x_{ijk}^{mno}\leq1,\forall i\in S,j\in G,\\
 k\in TW_{mn}^{ij},m\in A_i^s,o\in T_i
 \end{array}
\end{equation}
\begin{equation}
\begin{array}{l}
 \sum_{n\in A_j^g}\overline{x_{ijjk}^{mno}}\leq1,\forall i\in S,j,j^{\prime}\in G,\\
 k\in TW_{mn}^{ij},m\in A_i^s,o\in\overline{T_i}
 \end{array}
\end{equation}

$\bullet$ Each link construction task can be executed at most once.
\begin{equation}
\begin{array}{l}
 \sum_{k\in TW_{mn}^{ij}}\sum_{m\in A_i^s}\sum_{n\in A_j^g}x_{ijk}^{mno}\leq1,\\
 \forall i\in S,j\in G,o\in T_i
 \end{array}
\end{equation}
\begin{equation}
\begin{array}{l}
 \sum_{k\in TW_{mn}^{ij}}\sum_{m\in A_i^s}\sum_{n\in A_j^g}\overline{x_{ijj^{\prime}k}^{mno}}\leq1,\\
 \forall i\in S,j,j^{\prime}\in G,o\in\overline{T_i}
 \end{array}
\end{equation}

$\bullet$ The range of values for the decision variable.
\begin{equation}
	x_{ijk}^{mno}\in\{0,1\}
\end{equation}
\begin{equation}
	\overline{x_{ijj^{\prime}k}^{mno}}\in\{0,1\}
\end{equation}

\textcolor[rgb]{0,0,0}{\textbf{-Constraints under regular working mode:}}

$\bullet$ After completing a link construction task, the satellite must wait for a minimum interval before it can proceed with the next task.
\begin{equation}
\begin{array}{l}
et_o\cdot x_{ijk}^{mno}+\tau_{im}\leq st_{o^{\prime}}+M\cdot\left(1-x_{ij^{\prime}k^{\prime}}^{mn^{\prime}o^{\prime}}\right),\\
\forall i\in S,j,j^{\prime}\in G,k\in TW_{mn}^{ij},k^{\prime}\in TW_{mn^{\prime}}^{ij^{\prime}},\\
o,o^{\prime}\in T_i
\end{array}
\end{equation}

$\bullet$ After completion of a link construction task, the ground station must adhere to the minimum interval requirement between antenna task switching before proceeding with the next task.
\begin{equation}
\begin{array}{l}
et_o\cdot x_{ijk}^{mno}+\varphi_{ii^{\prime}}^j\leq st_{o^{\prime}}+M\cdot\left(1-x_{i^{\prime}jk^{\prime}}^{m^{\prime}no^{\prime}}\right),\\
\forall i\in S,j,j^{\prime}\in G,k\in TW_{mn}^{ij},k^{\prime}\in TW_{m^{\prime}n}^{i^{\prime}j},\\
o\in T_i,o^{\prime}\in T_{i^{\prime}}
\end{array}
\end{equation}

\textcolor[rgb]{0,0,0}{\textbf{-Constraints under feed switching operation mode:}}

$\bullet$ The duration of overlap between the two time windows for satellite feed switching exceeds the minimum required duration.
\begin{equation}
\begin{array}{l}
 \Delta_{innm^{\prime}}+M\cdot\left(2-\overline{x_{ijj^{\prime}k}^{mno}}-\overline{x_{ijj^{\prime}k^{\prime}}^{m^{\prime}n^{\prime}o^{\prime}}}\right)\\
 \leq\left|\overline{st_{o^{\prime}}}\cdot\overline{x_{ijj^{\prime}k^{\prime}}^{m^{\prime}n^{\prime}o^{\prime}}}-\overline{et_{o}}\cdot\overline{x_{ijj^{\prime}k}^{mno}}\right|,\\
 \forall i\in S,j,j'\in G,k\in TW_{mn}^{ij},\\
 k'\in TW_{m'n'}^{ij'},o,o'\in\overline{T_{i}}
 \end{array}
\end{equation}

$\bullet$ The temporal constraints of the task for constructing the power feeding switching link.
\begin{equation}
\begin{array}{l}
\overline{et_o}\cdot\overline{x_{ijj^{\prime}k}^{mno}}-\overline{st_{o^{\prime}}}\cdot\overline{x_{ijj^{\prime}k^{\prime}}}=0,\\
\forall i\in S,j,j^{\prime}\in G,k\in TW_{mn}^{ij},\\
k^{\prime}\in TW_{m^{\prime}n^{\prime}}^{ij^{\prime}},o,o^{\prime}\in\overline{T_i}
\end{array}
\end{equation}

\section{Proposed Method}
\label{Proposed Method}

To address the CS-GSNSP, we make corresponding enhancements to the regular genetic algorithm framework and propose an evolutionary algorithm (FFEEA) that incorporates fuzzy methodology for fitness function evaluation. The novel approach simplifies the original fitness evaluation process, with only a limited number of individuals using the actual fitness evaluation method within the model. Additionally, our algorithm employs adaptive population evolution operations to achieve a highly efficient search. This section provides an overview of the overall framework of FFEEA and outlines the improvement methods adopted in our algorithm.

\subsection{Overall Framework}

FFEEA, using fuzzy fitness evaluation, presents a novel approach to address diverse combinatorial optimization problems such as the CS-GSNSP. \textcolor[rgb]{0,0,0}{For problems such as CS-GSNSP, the algorithmic search often produces suboptimal solutions. Therefore, it is impractical to invest significant effort in evaluating these solutions. Instead, employing fuzzy fitness evaluation can significantly mitigate the cost associated with assessing low-quality solutions.} This algorithm builds upon the traditional genetic algorithm by incorporating further enhancements. Genetic algorithms are population-based evolutionary algorithms renowned for their robust exploration capabilities and effective handling of complex problems. These merits motivate us to propose a novel evolutionary method for solving CS-GSNSP. Our algorithm introduces three key innovations: 1) FFEM is used in the algorithm as an adaptive choice of fuzzy or real evaluation method based on individual performance (i.e., the total length of time of the chain-building task). Specifically, only individuals exhibiting high membership relationships with the optimal individual employ the model in Section \ref{Fuzzy Fitness Evaluation Method (FFEM)} to calculate actual fitness function values. 2) An adaptive crossover method is designed to automate the process of selecting from multiple operators to search for a better solution for task execution. 3) An elite retention strategy is introduced to accelerate convergence during the initial phase of algorithmic search. It can also ensure population diversity by discontinuing its usage after a certain number of iterations.

\subsection{Encoding and Decoding}

The coding method plays a crucial role in evolutionary algorithms, as it directly impacts individual fitness evaluation and scheme generation. In the FFEEA, we employ integer encoding to transform the link construction task sequence into individual chromosomes, where each gene corresponds to a specific task number. To facilitate comprehension, we present a simple example for intuitive illustration.

\begin{figure}[htp]
\centering
\includegraphics[width=0.45\textwidth]{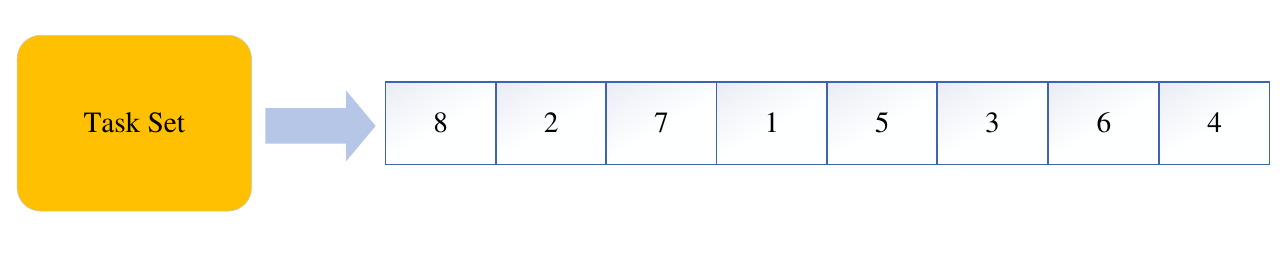}
\caption{An Example of Encoding}
\label{An Example of Encoding}
\end{figure}

\textcolor[rgb]{0,0,0}{As depicted in Figure \ref{An Example of Encoding}, a sequence of 8 link-building tasks is randomly generated to form an individual chromosome. Within this sequence, the gene coding value "8" represents the specific task related to the 8-th, while the coding value "2" corresponds to the 2-nd task. The remaining coding values remain unchanged.}

The process of evaluating fitness using either real models or fuzzy methods is commonly referred to as the decoding process, wherein the profit value is calculated based on the task execution scheme under constrained conditions. When employing a real model, the profit value derived from Eq. \ref{actual profit} is dependent on the task execution time. In this context, tasks are prioritized according to the maximum time allocation principle, whereby link construction tasks are executed for as long as feasible within the operational capacity of communication satellite and ground station equipment.

\subsection{Fuzzy Fitness Evaluation Method (FFEM)}
\label{Fuzzy Fitness Evaluation Method (FFEM)}
Complex combinatorial optimization problems often incur significant computational costs during fitness evaluation, particularly in evolutionary algorithms where each individual's fitness needs to be re-evaluated after structural improvements. During population evolution, numerous individuals may exhibit similar performance and struggle to surpass local optimal solutions, making it imprudent to allocate excessive computing resources to these individuals. To mitigate this issue, we propose a fuzzy fitness evaluation method (FFEM). The core concept of FFEM involves selecting a baseline individual for evaluation and determining the evaluation approach for other individuals based on their similarity relationships. The fitness function is computed using either fuzzy evaluation functions (i.e. membership functions) or real evaluation functions depending on individual similarity levels. A detailed description of FFEM is provided in the subsequent section.

Firstly, it is essential to identify a base individual as the central point, denoted as $C$. \textcolor[rgb]{0,0,0}{This baseline individual will serve as the basis for determining the specific method employed to calculate individual fitness values. Subsequently, the optimal individual in the population search is designated as the baseline individual $C$. Following this, we compute the similarity between the new individual $X_i$ and the baseline individual using the following equation (i.e. intuitive Gaussian membership function).}
\begin{equation}
	\label{similarity}
 \mu_i=\sum_{j\in T}\frac{\exp\left(-\left(x_j^i-c_j\right)^2/\sigma_j^2\right)}{\left|T\right|}
\end{equation}
where $\sigma_j$ indicates the control parameter. The specific calculation equation is:
\begin{equation}
\sigma_j=\gamma\frac1{\left(e^{f(c_j)}\right)^\tau}
\end{equation}
where $\tau>0$ is the emphasis operator and $\gamma$ is the constant. $\gamma$ is generally set to 1.

Subsequently, the fitness function is determined based on the similarity between individuals, as depicted in Eq. \ref{fuzzy eval}. In case the membership function value of the similarity degree equals 1, an accurate calculation of the fitness function is performed using the real evaluation function. Conversely, if the membership function value of similarity falls below 1, a rapid evaluation of fitness is conducted using both base individual and membership functions. The precise procedure for obtaining the fitness function value through fuzzy evaluation can be observed in Eq. \ref{fuzzy eval}.
\begin{equation}
\label{fuzzy eval}
f\left( {{X_i}} \right) = \left\{ {\begin{array}{*{20}{l}}
{\widehat f\left( {{X_i}} \right)}&{{\rm{if }} \ {\mu _i} < 1}\\
{f\left( {{X_i}} \right)}&{{\rm{otherwise}}}
\end{array}} \right.
\end{equation}
\textcolor[rgb]{0,0,0}{where $f\left( {{X_i}} \right)$ denotes the real evaluation function, $\hat{f}(X_i)$ refers to the following function.}
\begin{equation}
\label{fuzzy profit}
\hat{f}(X_i)=f(C)\cdot\mu_i
\end{equation}
where $\mu_i$ indicates the similarity, $f(C)$ indicates the fitness function value of the baseline individual.

If the fitness function value obtained through real fitness evaluation exceeds that of the base individual, both the base individual and its corresponding fitness function value will be updated. The FFEM simplifies the fitness evaluation process but may also introduce new challenges. Specifically, if this evaluation method is consistently applied, it could potentially overlook a portion of high-quality individuals. To mitigate such occurrences, an $\varepsilon$-evaluation strategy selection mechanism is employed. This mechanism compares random numbers $rand()$ and $\varepsilon$ within the range [0,1] to determine whether to directly calculate the fitness function value using real evaluation functions. Based on these explanations, Algorithm Table \ref{Fuzzy Fitness Evaluation Method (FFEM)} presents the pseudo-code for implementing FFEM.

\begin{algorithm}[htbp]
	\caption{Fuzzy Fitness Evaluation Method (FFEM)}
	\label{Fuzzy Fitness Evaluation Method (FFEM)}
	\LinesNumbered
	\KwIn{ population $P$, \textcolor[rgb]{0,0,0}{baseline individual} $C$, membership function $\widehat f\left( {{X_i}} \right)$,threshold $\varepsilon$}
	\KwOut{ fitness function value $F$}
	$Part_1,Part_2 \leftarrow$ Identify Individuals' evaluation methods based on the  comparison between $rand()$ and $\varepsilon$\;
	$Part'_1,Part_3 \leftarrow$ Divide offsprings in $Part_1$ according to similarity by Eq. \ref{similarity}\;
	$F_1 \leftarrow$ Fast estimate of the fitness function value in $Part'_1$ by Eq. \ref{fuzzy profit}\;
	$Part'_2 \leftarrow$ Merge $Part_2$ and $Part_3$\;
	$F_2 \leftarrow$ Evaluate fitness in $\left(Part'_2\right)$ using the actual objective function by Eq. \ref{actual profit}\;
	$F \leftarrow$ Merge the values of the fitness function of the two offspring populations$\left(F_1,F_2\right)$\;
\end{algorithm}

\subsection{Adaptive Crossover}

 In addition to employing fuzzy fitness evaluation, an adaptive crossover approach is also used in FFEEA. Crossover serves as a crucial operation for population evolution in genetic algorithm, aiming to discover high-quality solutions. This section initially presents various crossover operators employed in FFEEA and subsequently introduces the fundamental steps of the adaptive crossover method. \textcolor[rgb]{0,0,0}{Four types of crossover operators have been devised: double-point crossover (short segment), double-point crossover (long segment), recombination, and flip}. The specific procedures for these four crossover operators are described as follows.

\begin{figure}[htbp]
\centering
\includegraphics[width=0.45\textwidth]{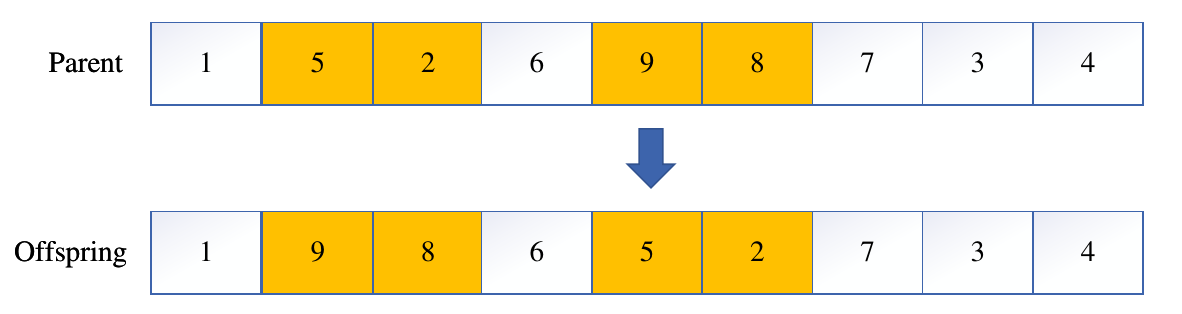}
\caption{An Example of Double-point Crossover (Short Fragment)}
\label{An Example of Double-point Crossover (Short Fragment)}
\end{figure}

\textbf{Double-point crossover (short fragment)}: Two gene fragments, each with a length of $L$, are randomly selected from an individual. Subsequently, these two fragments undergo a swapping operation to generate a novel individual (refer to Figure \ref{An Example of Double-point Crossover (Short Fragment)}).

\begin{figure}[htbp]
\centering
\includegraphics[width=0.45\textwidth]{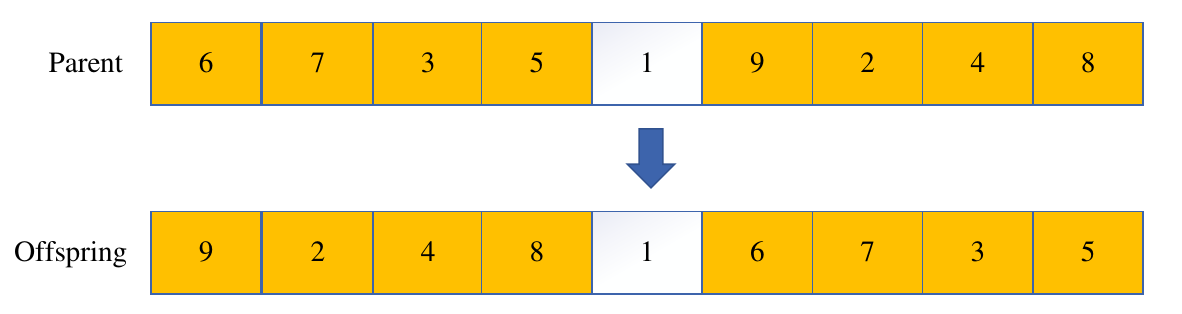}
\caption{An Example of Double-point Crossover (Long Fragment)}
\label{An Example of Double-point Crossover (Long Fragment)}
\end{figure}

\textbf{Double-point crossover (long fragment)}: Two gene fragments, each with a length of $2L$, are randomly selected from an individual. Subsequently, these two fragments undergo a swapping process to generate a novel individual (refer to Figure \ref{An Example of Double-point Crossover (Long Fragment)}).

\begin{figure}[htbp]
\centering
\includegraphics[width=0.45\textwidth]{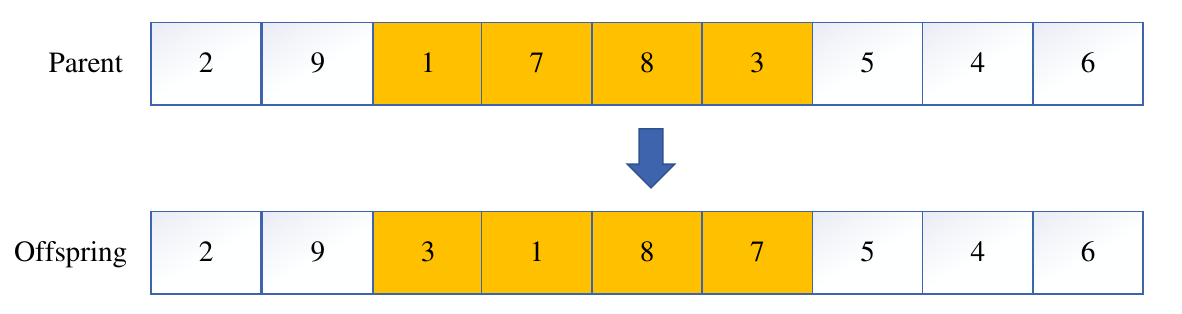}
\caption{An Example of Recombination}
\label{An Example of Recombination}
\end{figure}

\textbf{Recombination}: Genes are randomly selected from an individual chromosome fragment of length $2L$. Subsequently, these genes are reintegrated into the corresponding position of the individual in a relative order, following disruption of the original sequence (refer to Figure \ref{An Example of Recombination}).

\begin{figure}[htbp]
\centering
\includegraphics[width=0.45\textwidth]{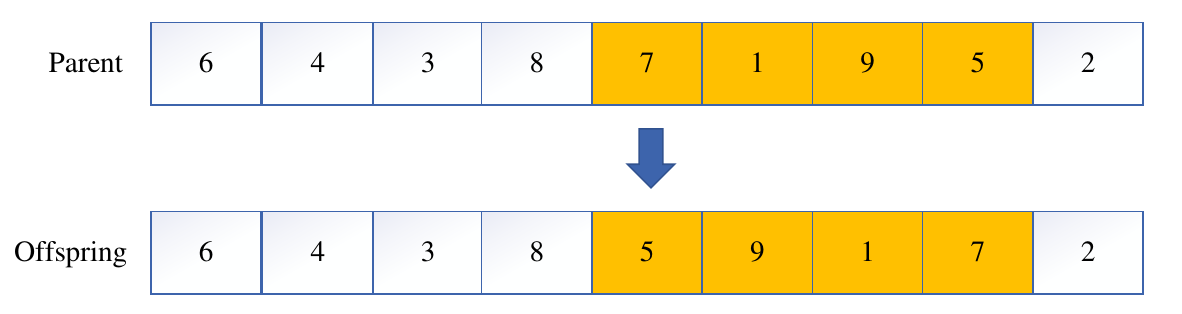}
\caption{An Example of Filp}
\label{An Example of Filp}
\end{figure}

\textbf{Flip}: A gene fragment of length $2L$ is randomly chosen from an individual, arranged in reverse order, and inserted back into the individual (refer to Figure \ref{An Example of Filp}).

After providing a comprehensive overview of the four types of crossover operators employed in FFEEA, this paper proceeds to elaborate on the subsequent steps including operator adaptive selection, score update, and weight update. Prior to commencing the population iterative search in FFEEA, all operators are initialized with identical scores, ensuring an equal probability of selection for each operator. During each round of iterative search, the roulette method is used for individual evolution mode selection. This approach enables every operator to be considered while favoring those with higher fitness values for individual evolution. The formula for operator adaptive selection is shown as follows.
\begin{equation}
\label{weight}
{w_i} = \frac{{so{c_i}}}{{\sum\nolimits_{i \in {O_c}} {so{c_i}} }}
\end{equation}
where ${w_i}$ indicates the weight of the $i$-th operator, ${so{c_i}}$ indicates the score of the $i$-th operator.

Then, the selected operator guides the evolutionary process within the contemporary population. Furthermore, based on fitness evaluation results from offspring populations, Algorithm Table \ref{Update Score} presents pseudo code illustrating how operator scores are updated.

\begin{algorithm}[htbp]
	\caption{Update Score}
	\label{Update Score}
	\LinesNumbered
	\KwIn{\textcolor[rgb]{0,0,0}{$last\_best$,$local\_best$,$global\_best$, selected operator score $sco_s$,acceptance percentage $\lambda$,$sco_1$, $sco_2$, $sco_3$}}
	\KwOut{ updated operator scores $sco_u$}
	\If{$local\_best>global\_best$}{
		$sco_u \leftarrow sco_s+sco_1$\;
	}
	\If{$local\_best \le global\_best$  and $local\_best>last\_best*\lambda$}{
		$sco_u \leftarrow sco_s+sco_2$\;
	}
	\If{$local\_best \le last\_best*\lambda$}{
		$sco_u \leftarrow sco_s+sco_3$\;
	}
\end{algorithm}

Subsequently, each specific evaluation times $\theta$ weight necessitates recalculation, with the dynamic adjustment operator being chosen accordingly. The Algorithm Table \ref{Calculate Weights of Operators} presents the pseudo-code for weight calculation. To mitigate excessive variation in individual weights, normalization of the operator weight is performed as follows.
\begin{equation}
	\label{normalize}
	sco_{i}^{'}=\frac{sco_{\max}-sco_{i}}{sco_{\max}-sco_{\min}}
\end{equation}
\begin{algorithm}[htbp]
	\caption{Calculate Weights of Operators}
	\label{Calculate Weights of Operators}
	\LinesNumbered
	\KwIn{crossover operator set $O_c$, score matrix \textbf{S}}
	\KwOut{ updated weight matrix \textbf{W}}
	$\textbf{S'} \leftarrow$Normalize scores by Eq. \ref{normalize}\;
	$sum \leftarrow$Calculate the sum of normalized scores$\left(\textbf{S'}\right)$\;
	\For{$i \leftarrow 1$ to $\left| {{O_c}} \right|$}{
		$w_i \leftarrow sco^n_i / sum$\;}
\end{algorithm}

\subsection{Elite Retention Strategy}

The algorithm should conduct an in-depth search within the search region surrounding individuals with high fitness values. To facilitate this, an elite retention strategy is devised to expedite the identification of optimal task execution. Specifically, at the onset of the search, if the population's best fitness value is lower than that discovered during the search process, FFEEA will replace the least performing individual with the best individual found since initiation. The elite retention strategy effectively guides the algorithm towards favorable outcomes but may also reduce population diversity. Consequently, once the threshold $Thre$ for fitness evaluations is reached, FFEEA will discontinue employing the elite retention policy.

\subsection{The FFEEA}
Building upon the foundation of traditional EA, our proposed FFEEA enhances the encoding and decoding process, fitness evaluation, and evolution operator based on the specific characteristics of the CS-GSNSP. This enhancement aims to optimize both search efficiency and algorithm performance. The pseudo-code for FFEEA is presented in Algorithm Table \ref{The FFEEA}.

\begin{algorithm}[htbp]
	\caption{The FFEEA}
	\label{The FFEEA}
	\LinesNumbered
	\KwIn{  $G$,$S$,$A$,$T$,$TW$,$N_p$,$Thre$,$\varepsilon$,$\lambda$,$\alpha$,$\beta$,$\theta$}
	\KwOut{ $Solution$}
	Initialize algorithm parameters\;
	$P_0 \leftarrow$ Generate initial population randomly\;
	$F,C \leftarrow$ Evaluate population fitness function values and generate center\;
	\While{not meet the termination condition}
	{
		$o_c \leftarrow$Select the operator adaptively by Eq. \ref{weight}\;
		$P' \leftarrow$ Perform individual selection and population crossover, mutation$\left(P\right)$\;
		$F \leftarrow$ Evaluate the value of the fitness function using FFEM in Algorithm \ref{Fuzzy Fitness Evaluation Method (FFEM)}\;
		$\textbf{S} \leftarrow$ Update scores using Algorithm \ref{Update Score}\;
		\eIf{$local\_best>global\_best$}{
			$global\_best \leftarrow local\_best$\;
			$global\_indi \leftarrow local\_indi$\;
			$C \leftarrow$Update center\;
		}
		{$count \leftarrow count+1$\;}
		$eavl \leftarrow$ Update the number of times the fitness function value was evaluated\;
            \If{$eval$ \% $\theta==0$ }{
            Calculate weights of operators using Algorithm \ref{Calculate Weights of Operators}\;
            }
		\If{$count<Thre$}{
			$P'' \leftarrow$ Use elite retention strategy$\left(P'\right)$\;
		}
            $last\_best \leftarrow local\_best$\;
	}
\end{algorithm}

\textcolor[rgb]{0,0,0}{As illustrated in Algorithm Table \ref{The FFEEA}, an initial population is generated randomly (Line 2). Subsequently, the algorithm employs roulette individual selection, adaptive crossover, and mutation to generate offspring populations. The fitness function values of the offspring population are evaluated using FFEM. In cases where a significant disparity exists between an individual and the optimal individual, it becomes necessary to use the mathematical model in Section \ref{Mathematical Model} for calculating the actual profit value achievable from the solution (Line 7). If an individual's actual profit value surpasses that of the optimal individual, both the optimal individual and its corresponding fitness value are updated (Line 10-12).}

\section{Experimental Studies}
\label{Experimental Studies}

This section presents the experimental validation of FFEEA's efficacy in addressing the CS-GSNSP, encompassing three components: experimental settings, experimental results, and discussion.

\begin{table*}[htbp]
	\caption{Max/Ave Results of Running Each Algorithm 30 Times}
	\label{Max/Ave Results of Running Each Algorithm 30 Times}
	\centering
	\begin{tabular}{lllllllll}
		\toprule[1.5pt]
		\multicolumn{1}{c}{\multirow{2}{*}{Instance}}& \multicolumn{2}{l}{FFEEA}         & \multicolumn{2}{l}{RL-GA}        & \multicolumn{2}{l}{KBGA}         & \multicolumn{2}{l}{ALNS-I}       \\ 
		\cmidrule(r){2-3} \cmidrule(r){4-5} \cmidrule(r){6-7} \cmidrule(r){8-9}
		\multicolumn{1}{c}{} & Max            & Ave              & Max           & Ave (WD)              & Max           & Ave (WD)               & Max           & Ave (WD)               \\
		\midrule[1pt]
		100-1    & \textbf{1295}  & \textbf{1295.0}  & \textbf{1295} & \textbf{1295.0=} & \textbf{1295} & \textbf{1295.0=} & \textbf{1295} & \textbf{1295.0=} \\
		100-2    & \textbf{1271}  & \textbf{1271.0}  & \textbf{1271} & \textbf{1271.0=} & \textbf{1271} & \textbf{1271.0=} & \textbf{1271} & \textbf{1271.0=} \\
		100-3    & \textbf{1231}  & \textbf{1231.0}  & \textbf{1231} & \textbf{1231.0=} & \textbf{1231} & \textbf{1231.0=} & \textbf{1231} & \textbf{1231.0=} \\
		200-1    & \textbf{2554}  & \textbf{2554.0}  & \textbf{2554} & \textbf{2554.0=} & \textbf{2554} & \textbf{2554.0=} & \textbf{2554} & \textbf{2554.0=} \\
		200-2    & \textbf{2471}  & \textbf{2471.0}  & \textbf{2471} & \textbf{2471.0=} & \textbf{2471} & \textbf{2471.0=} & \textbf{2471} & \textbf{2471.0=} \\
		200-3    & \textbf{2375}  & \textbf{2375.0}  & \textbf{2375} & \textbf{2375.0=} & \textbf{2375} & \textbf{2375.0=} & \textbf{2375} & \textbf{2375.0=} \\
		300-1    & \textbf{3822}  & \textbf{3821.1}  & 3819          & 3812.8-          & 3818          & 3817.0-          & 3815          & 3790.2-          \\
		300-2    & \textbf{3692}  & \textbf{3691.9}  & \textbf{3692} & 3686.5-          & \textbf{3692} & 3691.2-          & 3689          & 3670.9-          \\
		300-3    & \textbf{3706}  & \textbf{3703.2}  & 3704          & 3700.7-          & 3705          & 3701.9-          & 3701          & 3686.2-          \\
		400-1    & \textbf{5042}  & \textbf{5042.0}  & \textbf{5042} & 5041.2-          & \textbf{5042} & 5041.8-          & 5041          & 5040.7-          \\
		400-2    & \textbf{5190}  & \textbf{5174.6}  & 5178          & 5160.4-          & 5185          & 5165.7-          & 5072          & 5070.3-          \\
		400-3    & \textbf{4937}  & \textbf{4937.0}  & \textbf{4937} & \textbf{4937.0=} & \textbf{4937} & \textbf{4937.0=} & \textbf{4937} & \textbf{4937.0=} \\
		500-1    & \textbf{6143}  & \textbf{6123.3}  & 6122          & 6106.8-          & 6134          & 6112.3-          & 6111          & 6085.9-          \\
		500-2    & \textbf{6254}  & \textbf{6254.0}  & \textbf{6254} & 6253.4-          & 6253          & 6252.8-          & 6253          & 6243.1-          \\
		500-3    & \textbf{6173}  & \textbf{6172.9}  & 6171          & 6169.9-          & \textbf{6173} & 6172.8-          & 6164          & 6153.0-          \\
		600-1    & \textbf{7376}  & \textbf{7346.5}  & 7336          & 7300.6-          & 7333          & 7313.2-          & 7310          & 7269.6-          \\
		600-2    & \textbf{7229}  & \textbf{7195.1}  & 7153          & 7119.0-          & 7162          & 7129.3-          & 7124          & 7074.8-          \\
		600-3    & \textbf{7543}  & \textbf{7522.8}  & 7504          & 7467.9-          & 7509          & 7483.7-          & 7497          & 7431.5-          \\
		700-1    & \textbf{8410}  & \textbf{8398.5}  & 8400          & 8387.5-          & 8408          & 8393.8-          & 8388          & 8359.8-          \\
		700-2    & \textbf{8831}  & \textbf{8816.9}  & 8824          & 8805.1-          & 8822          & 8806.7-          & 8459          & 8403.1-          \\
		700-3    & \textbf{8565}  & \textbf{8565.0}  & \textbf{8565} & \textbf{8565.0=} & \textbf{8565} & \textbf{8565.0=} & 8563          & 8562.7-          \\
		800-1    & \textbf{9972}  & \textbf{9969.7}  & 9967          & 9959.8-          & 9968          & 9964.8-          & 9681          & 9928.4-          \\
		800-2    & \textbf{10080} & \textbf{10059.3} & 10051         & 10022.3-         & 10059         & 10026.4-         & 10020         & 9981.2-          \\
		800-3    & \textbf{10045} & \textbf{10043.7} & 10038         & 10032.9-         & 10041         & 10039.8-         & 10037         & 10008.2-         \\
		900-1    & \textbf{11637} & \textbf{11611.2} & 11562         & 11533.3-         & 11584         & 11547.5-         & 11547         & 11475.8-         \\
		900-2    & \textbf{11126} & \textbf{11100.2} & 11053         & 11012.5-         & 11076         & 11029.2-         & 11033         & 10956.5-         \\
		900-3    & \textbf{11582} & \textbf{11547.7} & 11510         & 11453.6-         & 11507         & 11466.9-         & 11461         & 11393.8-         \\
		1000-1   & \textbf{12669} & \textbf{12602.4} & 12503         & 12403.5-         & 12471         & 12421.8-         & 12502         & 12335.2-         \\
		1000-2   & \textbf{12458} & \textbf{12398.5} & 12328         & 12258.6-         & 12371         & 12287.0-           & 12270         & 12191.1-         \\
		1000-3   & \textbf{12181} & \textbf{12128.1} & 12011         & 11963.4-         & 12044         & 11983.8-         & 11956         & 11896.4-         \\
		\midrule[1pt]
		\multicolumn{3}{c}{+/-/=}                    & \multicolumn{2}{c}{0/22/8}       & \multicolumn{2}{c}{0/22/8}       & \multicolumn{2}{c}{0/23/7}  \\
		\bottomrule[1.5pt]
	\end{tabular}
\end{table*}

\subsection{Experimental Setting}
Experimental environment: All experiments are conducted with Core i7-1260P 2.10GHz CPU, 16GB memory capacity, and Matlab 2023a coding environment.

Test instances: Since there is no public test set for CS-GSNSP, we use a random method to generate test instances. The relevant parameters in the scenarios are set as follows: the time range of the scenarios is 24 hours, and the number of tasks is between 100-1000. The link construction task requirements time follows a normal distribution with a mean of 50 and a standard deviation of 40. The task profits follow a uniform distribution with a mean of 15 and a standard deviation of 8. The antenna attitude adjustment time is 10s, minimum overlap time for performing feed switching is 8s.

Comparison algorithms: \textcolor[rgb]{0,0,0}{In order to validate the performance of the proposed algorithm in addressing the CS-GSNSP, we employ a range of comparison algorithms that have been previously introduced in recent literature on satellite scheduling.} These comparison algorithms include a reinforcement learning-based genetic algorithm (RL-GA) \cite{song2023rl}, a knowledge-driven genetic algorithm (KBGA) \cite{song2020knowledge}, and an improved adaptive large neighborhood search algorithm (ALNS-I) \cite{liu2024knowledge}. RL-GA dynamically selects crossover operators using reinforcement learning techniques, KBGA uses domain knowledge-driven algorithms for efficient exploration, while ALNS-I incorporates specialized damage and repair operators tailored for satellite scheduling problems.

Algorithm parameter settings: \textcolor[rgb]{0,0,0}{For FFEEA, we systematically select various parameter configurations for the experiment based on prior research studies. The most stable combination of parameters is used as the algorithm's setting in the experimental setup.} In FFEEA, $N_p$ is set to 10, the maximum number of evaluations is set to 5000, $\alpha$ is set to 0.9, $\beta$ is set to 0.1, $\varepsilon$ is set to 0.9, $Thre$ is set to 500, $\theta$ is set to 100, $sco_1$ is set to 30, $sco_2$ is set to 20, $sco_3$ is set to 10, $\lambda$ is set to 0.95, $\gamma$ is set to 1, $\tau$ is set to 0.05. The parameter settings of comparison algorithms are consistent with those in the relevant literature.

Evaluation index: The experiment is conducted with each algorithm being run 30 times, and the recorded results include the maximum value (denoted as Max), average value (denoted as Ave), and standard deviation (denoted as Std) of a single run. Furthermore, a rank sum test (denoted as WD) at a significance level of $p=0.05$ is employed to examine potential significant differences among the results.

\subsection{Experimental Results}

\begin{figure*}[htbp]
	\centering
	\subfigure[Average CPU Time for Instance 400-2]{\includegraphics[width=0.45\textwidth]{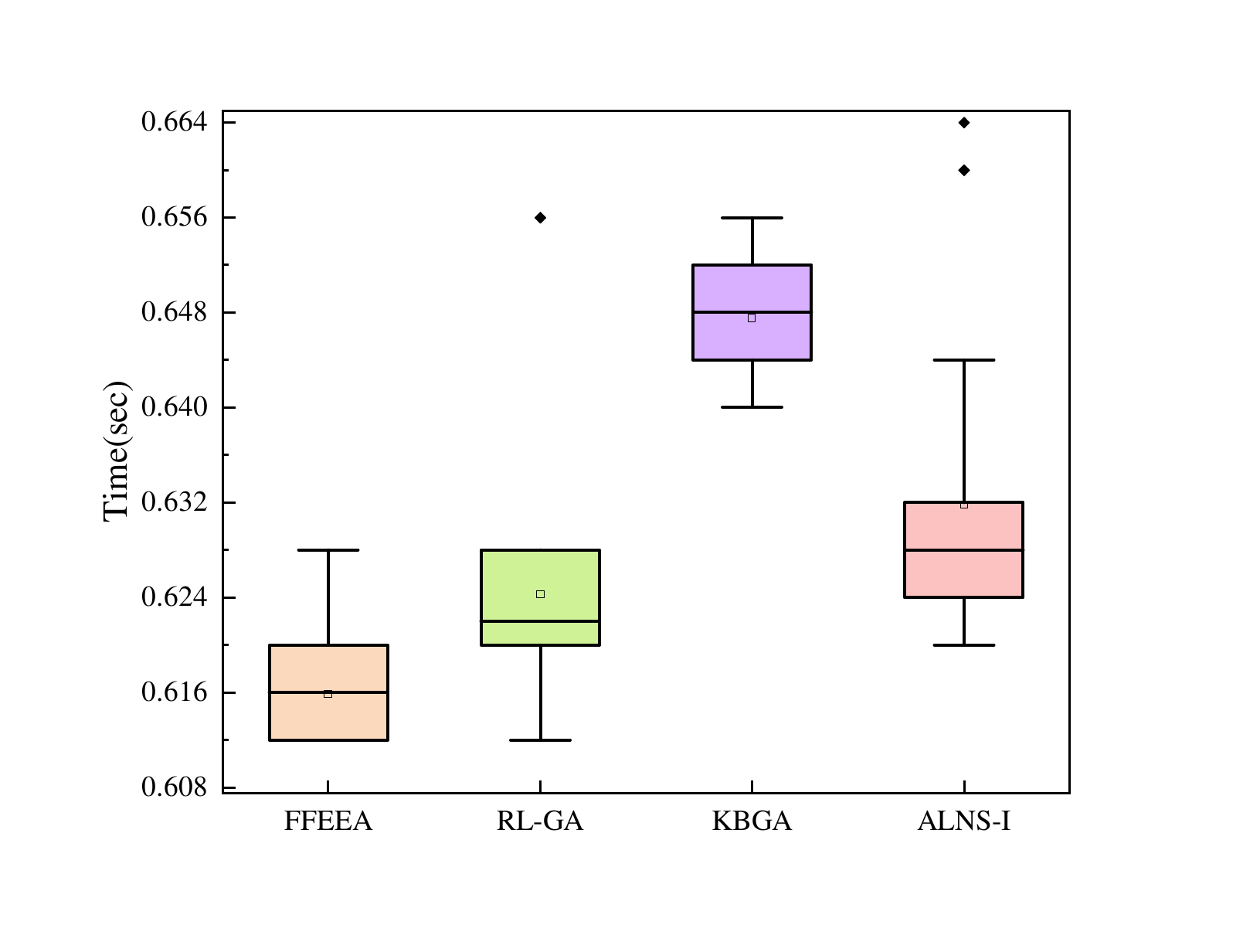}} \qquad
	\subfigure[Average CPU Time for Instance 800-2]{\includegraphics[width=0.45\textwidth]{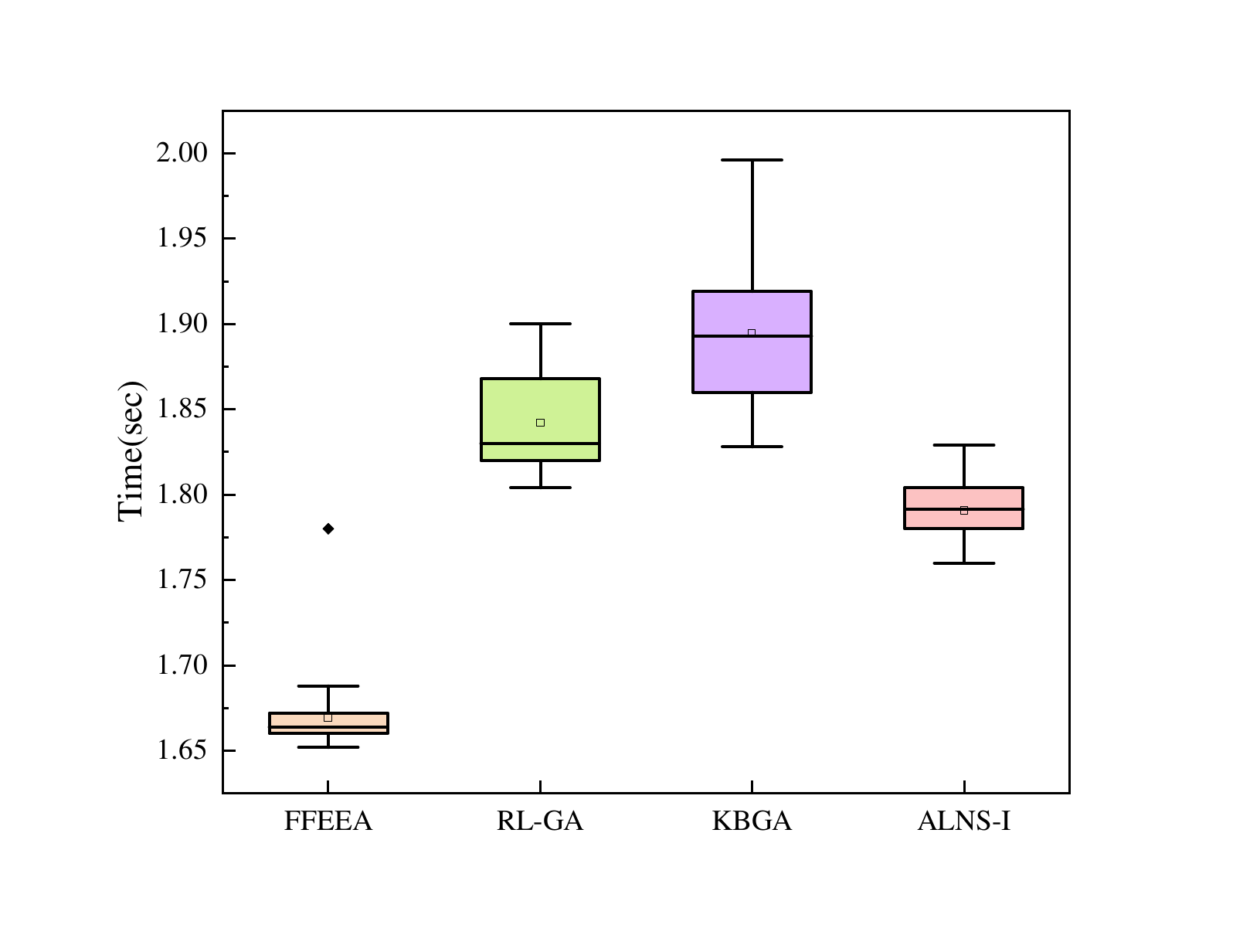}}
	\caption{Boxplot of Average CPU Time}
	\label{Boxplot of Average CPU Time}
\end{figure*}

\begin{figure*}[htbp]
	\centering
	\subfigure[Comparison on the Stability of 600 Task Scale]{
		\includegraphics[width=0.45\textwidth]{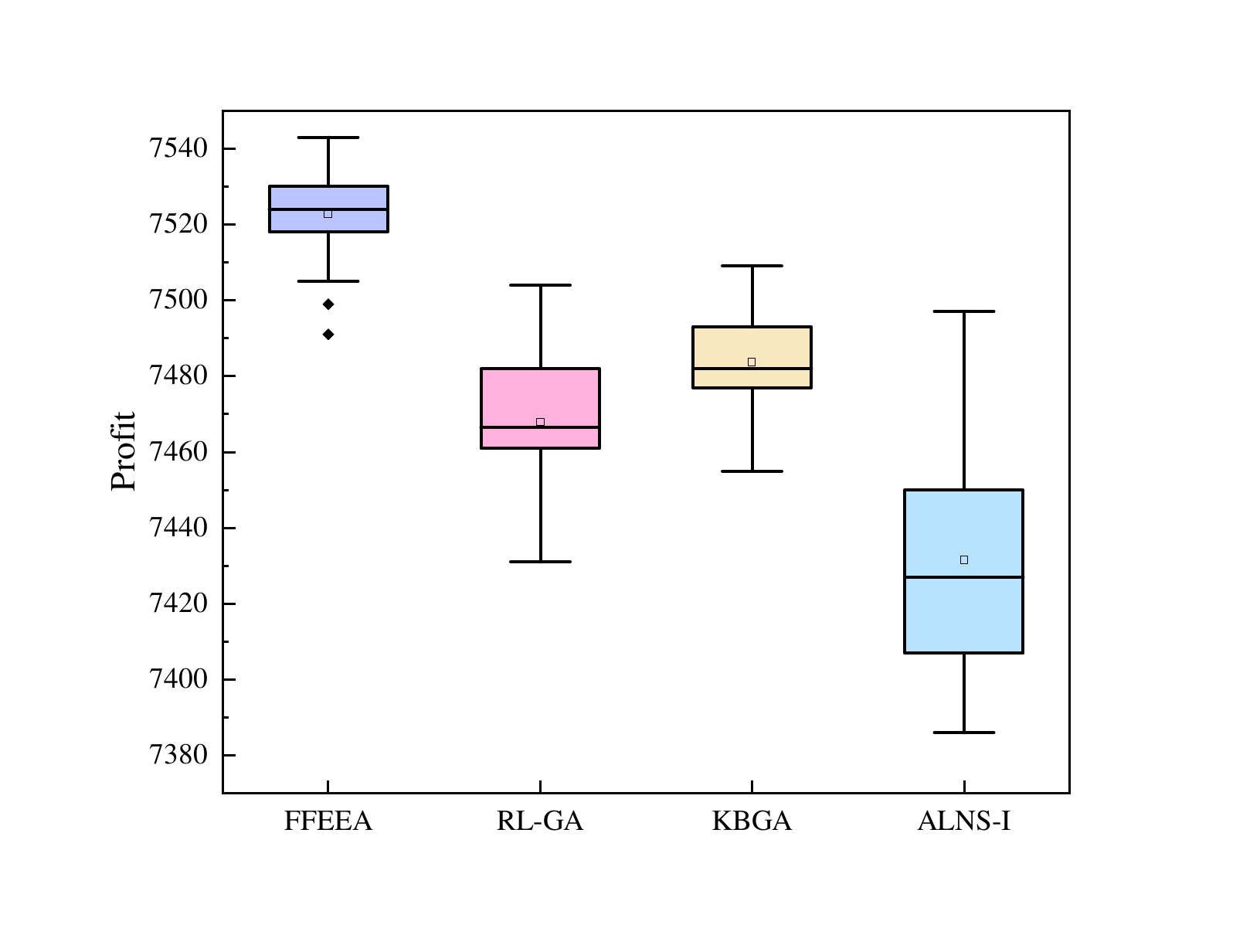}}
	\subfigure[Comparison on the Stability of 900 Task Scale]{
		\includegraphics[width=0.45\textwidth]{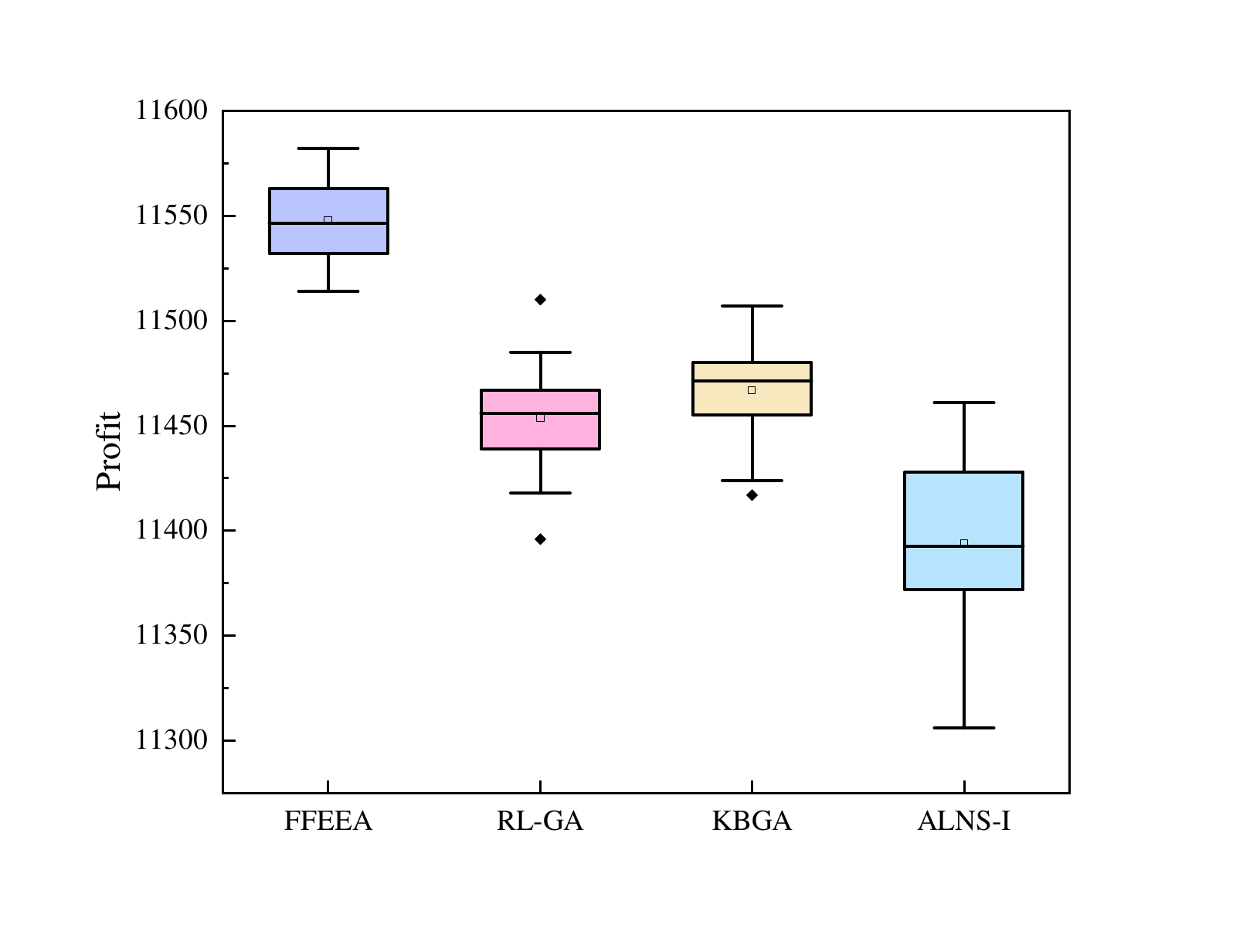}}
	\caption{Comparison on the Stability of Algorithms}
	\label{Comparison on the Stability of Algorithms}
\end{figure*}

\begin{figure*}[htbp]
	\centering
	\subfigure[Convergence Curves of 800 Task Scales]{
		\includegraphics[width=0.45\textwidth]{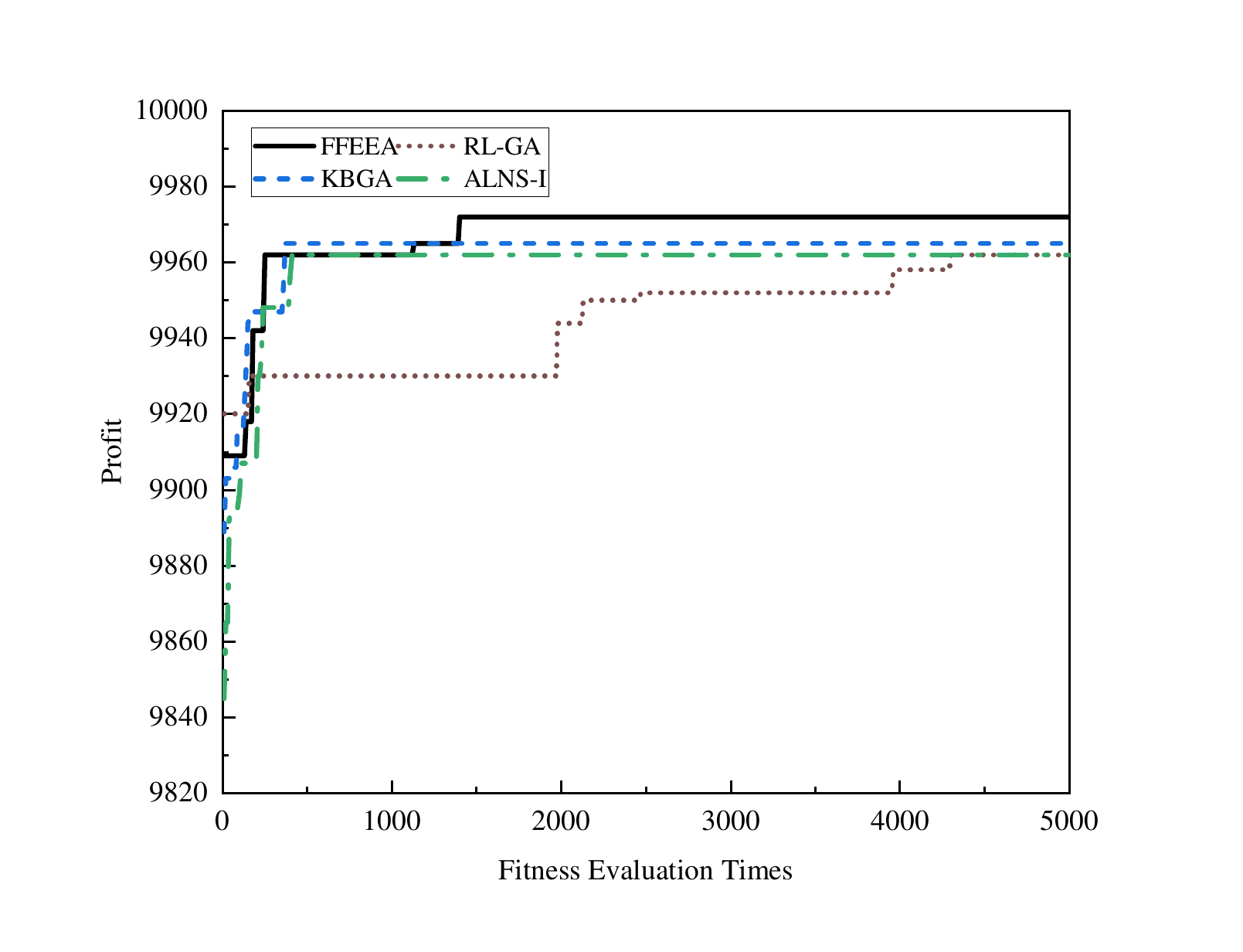}}
	\subfigure[Convergence Curves of 1000 Task Scales]{
		\includegraphics[width=0.45\textwidth]{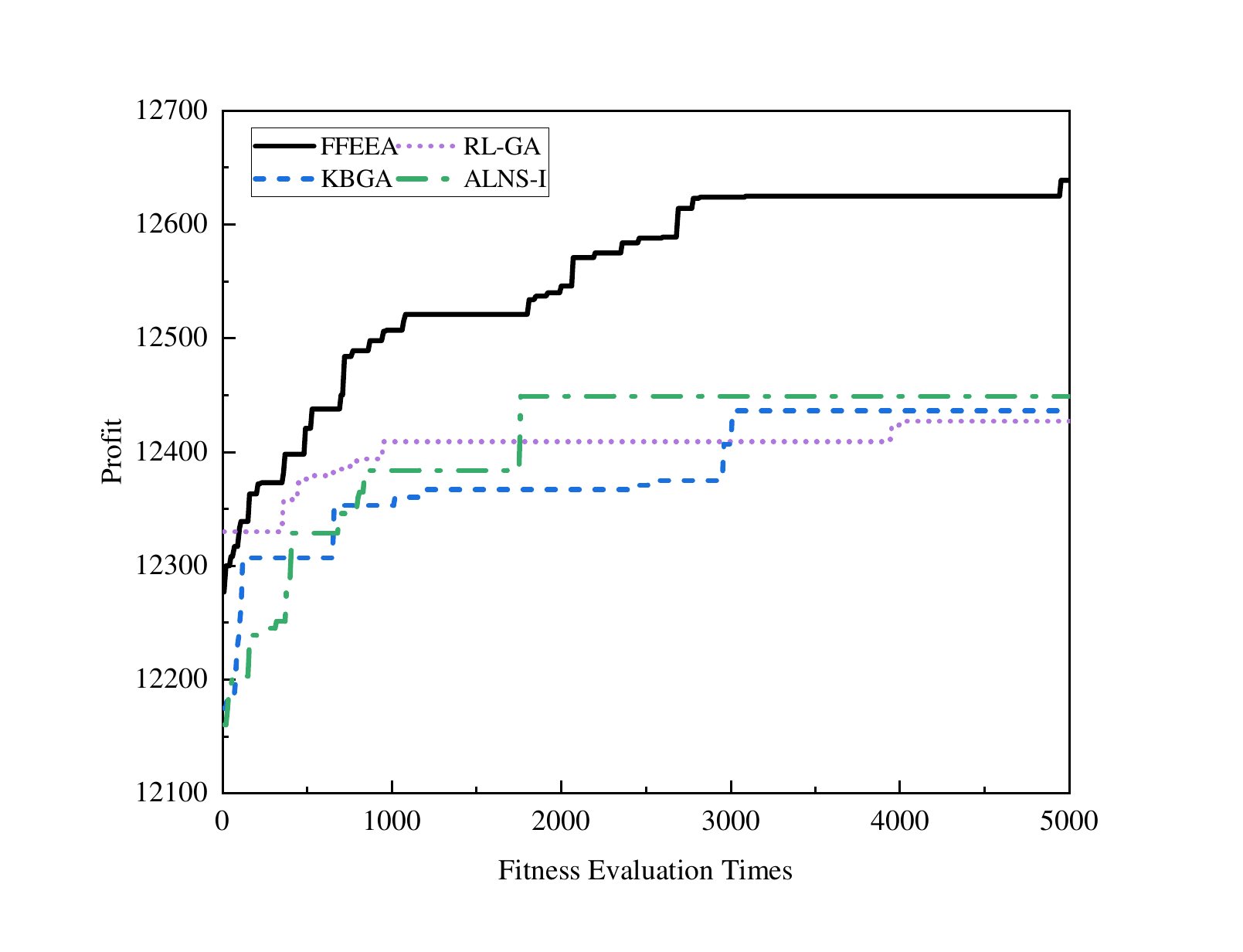}}
	\caption{Convergence Curves in 800 and 1000 Task Scales of Instances}
	\label{Convergence Curves in 800 and 1000 Task Scales of Instances}
\end{figure*}

\begin{figure*}[htbp]
	\centering
	\subfigure[Results of 800 Task Scales]{
		\includegraphics[width=0.45\textwidth]{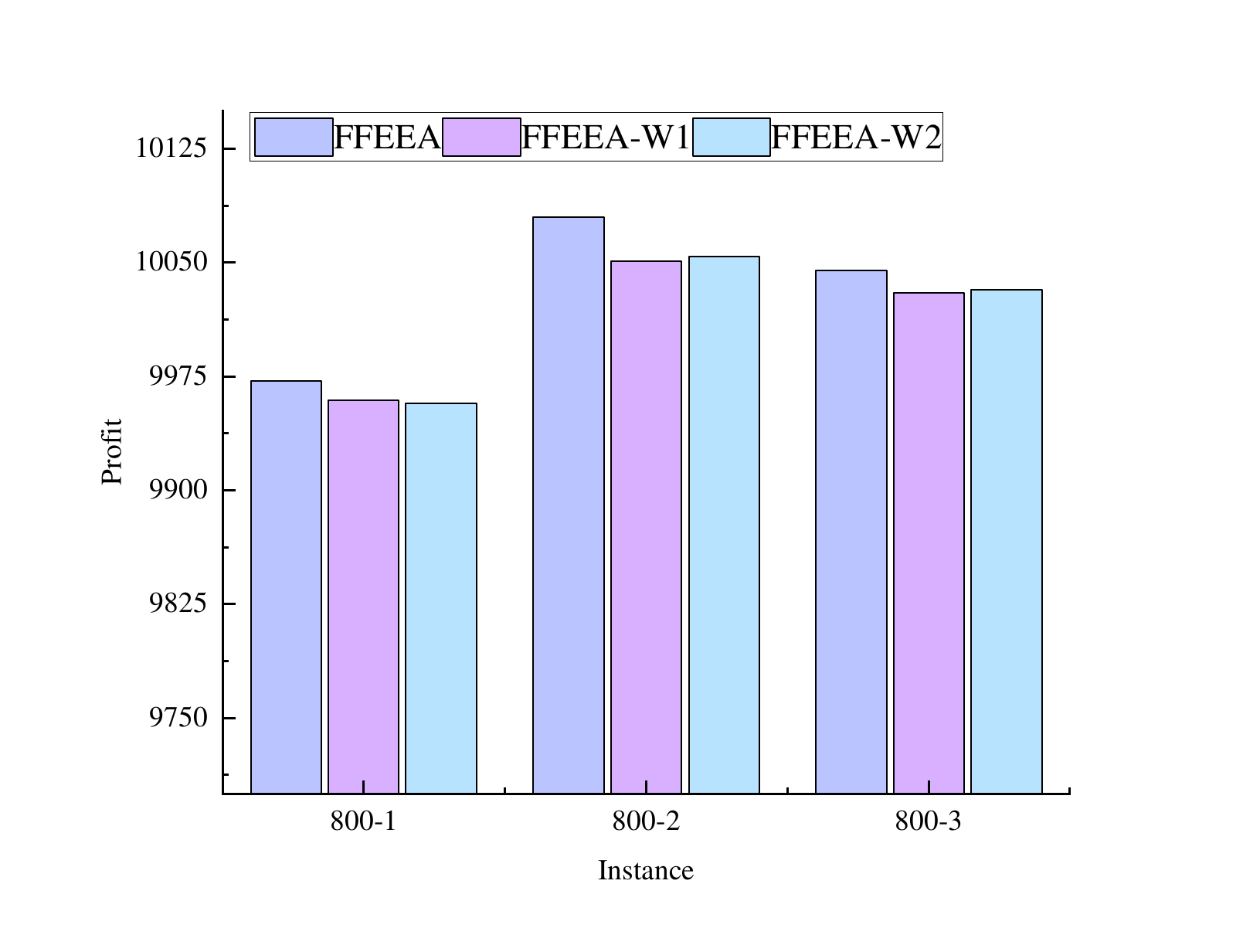}}
	\subfigure[Results of 1000 Task Scales]{
		\includegraphics[width=0.45\textwidth]{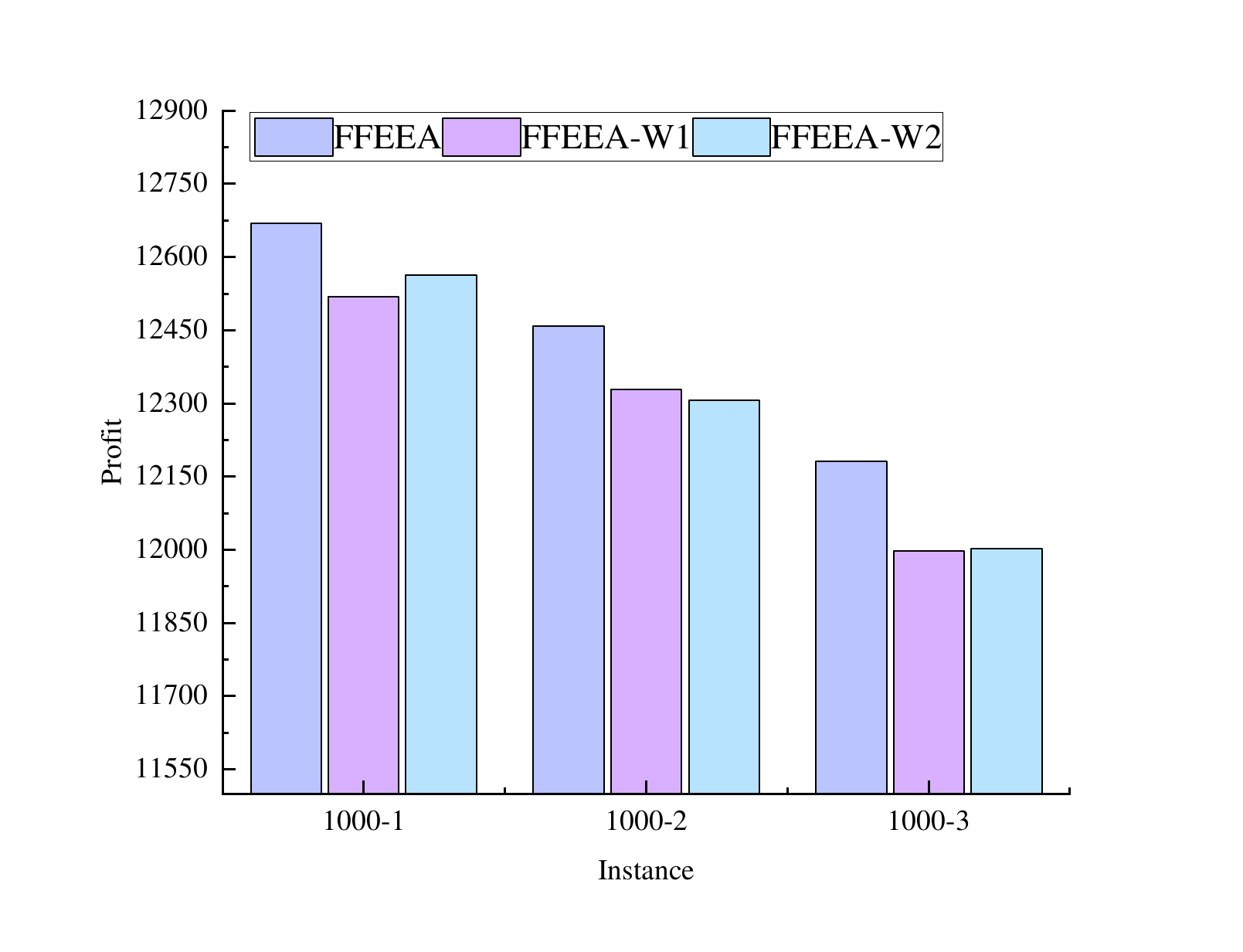}}
	\caption{Comparison Results of Algorithms with Different Strategies}
	\label{Comparison Results of Algorithms with Different Strategies}
\end{figure*}

Firstly, we compare the search efficiency of algorithms for different task sizes. As depicted in Table \ref{Results of Average CPU Time}, FFEEA demonstrates a shorter completion time for the set number of experiments. This indicates that the use of FFEM in the algorithm effectively reduces the time required for evaluating individual performance. Other algorithms employing real fitness evaluation methods exhibit similar timings, highlighting that fitness evaluation is a significant computational resource component for combinatorial optimization problems such as CS-GSNSP. The boxplot \ref{Boxplot of Average CPU Time} illustrates the CPU time of various algorithms applied to the instances 400-2 and 800-2. From this figure, it can be seen that FFEEA's execution times are relatively close across multiple runs, while other algorithms display noticeable variations due to distinct processing constraints during fitness evaluation.

Subsequently, the optimization performance of the algorithm is presented in Table \ref{Max/Ave Results of Running Each Algorithm 30 Times}. It can be seen that the FFEEA outperforms comparison algorithms in terms of both best and average values. The FFEM component effectively avoids neglecting high-quality individuals by simplifying fitness evaluation for certain individuals. Overall, as instance scale increases, the advantages of FFEEA over other algorithms become increasingly apparent. Among various comparison algorithms, RL-GA ranks highest followed by KBGA and ALNS-I.

\begin{table}[htbp]
\caption{Results of Average CPU Time}
\label{Results of Average CPU Time}
\centering
\begin{tabular}{lllllllll}
\toprule[1.5pt]
Instance & FFEEA          & RL-GA          & KBGA           & ALNS-I         \\ 
\midrule[1pt]
100-1    & \textbf{0.113} & 0.119          & 0.143          & 0.124          \\
100-2    & \textbf{0.111} & 0.114          & 0.140          & 0.119          \\
100-3    & \textbf{0.107} & 0.111          & 0.126          & 0.112          \\
200-1    & \textbf{0.209} & 0.212          & 0.233          & 0.220          \\
200-2    & \textbf{0.221} & 0.228          & 0.243          & 0.227          \\
200-3    & \textbf{0.218} & 0.219          & 0.239          & 0.224          \\
300-1    & \textbf{0.347} & 0.355          & 0.372          & 0.352          \\
300-2    & \textbf{0.363} & 0.369          & 0.391          & 0.375          \\
300-3    & 0.357          & \textbf{0.356} & 0.386          & 0.364          \\
400-1    & \textbf{0.606} & 0.625          & 0.975          & 0.983          \\
400-2    & \textbf{0.616} & 0.624          & 0.648          & 0.632          \\
400-3    & \textbf{0.891} & 0.895          & 1.020          & 0.952          \\
500-1    & \textbf{1.110} & 1.114          & 1.263          & 1.222          \\
500-2    & \textbf{1.086} & 1.098          & 1.202          & 1.387          \\
500-3    & \textbf{0.772} & 0.785          & 0.816          & 0.792          \\
600-1    & \textbf{1.021} & 1.022          & 1.037          & 1.038          \\
600-2    & 1.051          & 1.053          & 1.058          & \textbf{1.047} \\
600-3    & \textbf{0.997} & 1.001          & 1.016          & 0.999          \\
700-1    & \textbf{1.398} & 1.403          & 1.415          & 1.404          \\
700-2    & \textbf{1.397} & 1.408          & 1.420          & 1.405          \\
700-3    & \textbf{1.430} & 1.434          & 1.445          & 1.434          \\
800-1    & \textbf{1.746} & 1.757          & 1.755          & 1.752          \\
800-2    & \textbf{1.669} & 1.842          & 1.894          & 1.791          \\
800-3    & \textbf{1.713} & 1.718          & 1.722          & 1.718          \\
900-1    & \textbf{1.986} & 1.992          & 1.999          & 1.993          \\
900-2    & \textbf{2.039} & 2.052          & 2.050          & 2.045          \\
900-3    & \textbf{2.013} & 2.028          & 2.020          & 2.021          \\
1000-1   & 2.424          & 2.438          & \textbf{2.423} & 2.430          \\
1000-2   & \textbf{2.416} & 2.436          & 2.417          & 2.419          \\
1000-3   & 2.398          & 2.403          & \textbf{2.394} & 2.402          \\ 
\bottomrule[1.5pt]
\end{tabular}
\end{table}

Then, the stability and convergence of the algorithm are subsequently analyzed. The boxplot \ref{Comparison on the Stability of Algorithms} illustrates the fluctuation in search performance of the algorithm across calculation instances at different scales, revealing that search of FFEEA is more likely to discover task execution schemes with comparable performance compared to comparison algorithms. Figure \ref{Convergence Curves in 800 and 1000 Task Scales of Instances} demonstrates the convergence performance of each algorithm as well. The proposed algorithm effectively enhances scheme profitability through swift improvements during the search process, facilitated by an adaptive crossover population evolution strategy that can be dynamically adjusted.

We conduct experiments to validate the efficacy of the improvement strategies in FFEEA. The two algorithms, FFEEA without adaptive crossover (denoted as FFEEA-W1) and FFEEA without FFEM (denoted as FFEEA-W2) are used for comparison. To ensure fairness, we set a predetermined search time as a termination condition for both comparison algorithms to produce results. \textcolor[rgb]{0,0,0}{As depicted in Figure \ref{Comparison Results of Algorithms with Different Strategies}, the employed improvement strategies in FFEEA significantly enhance the algorithm's capability to explore high-quality solutions.} The roles of adaptive crossover and FFEM exhibit relatively similar effects within the algorithm.

\subsection{Discussion}

In our study, the proposed FFEEA offers an efficient approach to address practical engineering problems, including the CS-GSNSP. When dealing with complex problems and large solution spaces, the fuzzy fitness evaluation method employed in the algorithm significantly reduces computation time while still generating high-quality solutions. Moreover, adaptive crossover enhances the algorithm's search performance. FFEEA integrates various operators and independently selects those with superior performance to efficiently explore the solution space. Each operator has a chance of being selected, thereby increasing search strategy diversity and facilitating both local and global exploration. This crossover method shares similarities with the fitness evaluation technique and exhibits low time complexity.

The use of this algorithm design approach will inevitably have a certain impact on the search performance. This is due to the omission of fitness calculations for some high-quality individuals in the FFEM method. To mitigate this issue, we have devised the $\varepsilon$-evaluation strategy selection mechanism, which allows these individuals to be directly used for calculating the fitness function value of the model. Based on our experimental results, it can be observed that the algorithm's fast search does not significantly compromise the quality of scheduling schemes. While we have achieved success thus far, there are still numerous areas that warrant further improvement. For instance, additional refinement can be made to enhance the fuzzy fitness evaluation mechanism. Furthermore, we can explore leveraging insights gained during problem-solving processes to determine specific FE methods and introduce auxiliary decision-making approaches.

\section{Conclusion}
\label{Conclusion}

The construction scheme of the communication link between the satellite and ground station in a communication satellite network is crucial for achieving data transmission and providing diverse services to users. While regional satellite network communication scheduling has received significant attention, the issue of CS-GSNSP remains largely unexplored. To enhance the service quality of CSN, we propose a mathematical programming model for addressing this problem. Subsequently, considering the extensive solution space and high search complexity associated with CS-GSNSP, we introduce an FFEM approach along with the design of FFEEA. The FFEM method adapts by selecting fitness evaluation strategies based on individual performance. Finally, numerous simulation experiments validate the effectiveness of our proposed algorithm in solving CS-GSNSP.

In future studies, we will conduct in-depth investigations into CS-GSNSP in specific cases, such as equipment failures of satellites or ground stations and the need for time-constrained communication link establishment. These complex practical scenarios impose higher demands on algorithm design. Additionally, it is essential to consider uncertainty factors in problem research, where distribution estimation methods and machine learning techniques serve as effective choices for problem-solving. Furthermore, we will further explore satellite-ground station data transmission scheduling by incorporating the satellite-ground link construction scheme as input. Moreover, other strategies like small population size and heuristic initialization methods can be employed to guide the search process. The genetic algorithm evolution process can also benefit from additional mechanisms such as reinforcement learning methods and integrated strategies based on voting mechanisms.

\section*{Acknowledgements}
This work was supported by the National Natural Science Foundation of China (723B2002).

\section*{Declaration of Competing Interest}

The authors declare that they have no known competing financial interests or personal relationships that could have appeared to
influence the work reported in this paper.

\bibliographystyle{unsrt}

\bibliography{mybib}



\end{document}